\documentclass[a4paper]{amsart}
\title[Hyperbolic Schwarz map]{
Hyperbolic Schwarz map for the hypergeometric differential equation}
\date{February 28, 2007}
\usepackage[dvips]{graphicx} 
\usepackage{amssymb}
\usepackage{rsfs}
\usepackage[usenames]{color}

\usepackage{amsthm}
\theoremstyle{plain}
 \newtheorem{theorem}{Theorem}[section]
 
 \newtheorem{lemma}[theorem]{Lemma}
 
 \theoremstyle{remark}
 
\numberwithin{equation}{section}
\numberwithin{figure}{section}

\newcommand{\Z}{\boldsymbol{Z}}
\newcommand{\R}{\boldsymbol{R}}
\newcommand{\C}{\boldsymbol{C}}
\newcommand{\PP}{\boldsymbol{P}}
\newcommand{\Her}{\operatorname{Her}}
\newcommand{\HH}{\boldsymbol{H}}
\newcommand{\Mon}{\operatorname{Mon}}
\def\transpose#1{\mathord{\mathopen{{\vphantom{#1}}^t}#1}} 
\author{Takeshi Sasaki}
\address[Sasaki]{%
   Department of Mathematics,
   Kobe University,
   Kobe 657-8501, Japan%
}
\email{sasaki@math.kobe-u.ac.jp}

\author{Kotaro Yamada}
\address[Yamada]{%
   Faculty of Mathematics,
   Kyushu University,
   Fukuoka 812-8581, Japan%
}
\email{kotaro@math.kyushu-u.ac.jp}

\author{Masaaki Yoshida}
\address[Yoshida]{%
   Faculty of Mathematics,
   Kyushu University,
   Fukuoka 810-8560, Japan%
}
\email{kotaro@math.kyushu-u.ac.jp}
\keywords{hypergeometric differential equation, Schwarz
          map, hyperbolic Schwarz map, flat surfaces, 
	  flat fronts}
\subjclass[2000]{33C05, 53C42}

\begin{document}
\maketitle
\begin{abstract}
 The Schwarz map of the hypergeometric differential equation is studied
 since the beginning of the last century. Its target is the complex
 projective line, the 2-sphere. This paper introduces the {\it
 hyperbolic Schwarz map}, whose target is the hyperbolic 3-space. This
 map can be considered to be a lifting to the 3-space of the Schwarz map. 
 This paper studies the singularities of this map, 
 and visualize its image when the monodromy group is a finite group 
 or a typical Fuchsian group. 
 General cases will be treated in forthcoming papers.
\end{abstract}

\tableofcontents

\section{Introduction}
Consider the {\it hypergeometric differential equation} 
\begin{equation}\tag*{$E(a,b,c)$}\label{eq:hgde}
   x(1-x)u''+\{c-(a+b+1)x\}u'-abu=0,
\end{equation}
and define its {\it Schwarz map} by
\begin{equation}\label{eq:schwarz}
   s:X=\C-\{0,1\}\ni x\longmapsto u_0(x):u_1(x)\in Z\cong\PP^1,
\end{equation}
where $u_0$ and $u_1$ are linearly independent solutions of
\ref{eq:hgde}
and $\PP^1$ is the complex projective line.
The Schwarz map of the hypergeometric differential equation was studied
by Schwarz when the parameters $(a,b,c)$ are real.

Since the work on Schwarz' study was a great success, 
many kinds of high-dimensional versions are
studied analytically, algebro-geometrically and arithmetically, in these
decades. But the authors have had a slight reservation about 
the Schwarz map in \eqref{eq:schwarz}: Its
target seems not to be exactly the correct one, 
because even if
the monodromy group of $s$, the projective monodromy group of the equation, 
is discrete in $PGL_2(\C)$, it does not,
in general, act properly discontinuously on any non-empty open set of
the target $\PP^1$, and so the image would be chaotic.

We propose a variation of the Schwarz map ({\it hyperbolic Schwarz map}), 
which solves the difficluty above. It is defined as follows: Change the
equation \ref{eq:hgde} into the so-called $SL$-form: 
\begin{equation} 
\tag{ESL}\label{eq:sl-form}
    u''-q(x)u=0,
\end{equation}
and transform it to the matrix equation
\begin{equation}\label{eq:sl-form-mat}
   \dfrac{d}{dx}(u,u')=(u,u')\Omega,\quad \Omega
        =\left(\begin{array}{cc}0&q(x)\\1&0\end{array}\right).
\end{equation}
We now define the {\it hyperbolic Schwarz map}, 
denoted by $\mathscr{S}$, as the composition of the (multi-valued) map
\begin{equation}\label{eq:hsm1}
    X\ni x\longmapsto H=U(x)\transpose{\overline U}(x)\in \Her^+(2)
\end{equation}
and the natural projection $\Her^+(2)\to \HH^3:=\mbox{Her}^+(2)/\R^+,$
where $U(x)$ is a fundamental solution of the system, 
$\Her^+(2)$  the space of positive-definite Hermitian matrices 
of size $2$, and $\R^+$ the multiplicative group of positive real numbers; 
the space $\HH^3$ is called the {\em hyperbolic $3$-space}. 
Note that the target of the hyperbolic Schwarz map is 
$\HH^3$, whose boundary is $\PP^1$, 
which is the target of the Schwarz map. In this sense, our hyperbolic
Schwarz map is a lift-to-the-air of the Schwarz map. 
Note also that the monodromy group of the system acts naturally on $\HH^3$.

Here we state a defect of our hyperbolic Schwarz map.
There are no standard way to transform our equation \ref{eq:hgde} into a
matrix system (this freedom is often called the {\it gauge ambiguity}); 
we made a choice.  The cost is that the symmetry of 
the equation \ref{eq:hgde}, 
which descends to the Schwarz map, does not necessarily descends to the
hyperbolic Schwarz map.

But, on the contrary, thanks to this choice, 
the image surface (of $X$ under $\mathscr{S}$) has 
the following geometrically nice property: It  is one of the 
{\it flat fronts} in $\HH^3$,
which is a flat surface with a certain kind of singularities
\cite{KUY2}.
Moreover, 
the classical Schwarz map $s$ is recovered  as the {\em hyperbolic 
Gauss map\/} of the hyperbolic Schwarz map as a flat front.
The papers \cite{GMM, KUY} gave a method
of constructing flat surfaces in the three-dimesnional hyperbolic space.
Since any closed nonsigular flat surface is isometric 
to a horosphere or a hyperbolic cylinder, such surfaces have necessarily 
singularities: 
generic singularities of flat fronts
are cuspidal edges and swallowtail singularities
\cite{KRSUY}, see Section~\ref{sec:sing}.

We will publish a series of papers about the hyperbolic 
Schwarz map and its singularities (\cite{SYY2, NSYY}). 
This is the first one. 
In this paper, we study the hyperbolic Schwarz map $\mathscr{S}$ 
of the equation
\ref{eq:hgde} when the parameters $(a,b,c)$ are real, especially when
its monodromy group is a finite (polyhedral) group or a Fuchsian group.
In general, generic singularities of flat fronts are 
cuspidal edges and swallowtails, 
In our special cases,
we find that, in each case, there is a simple closed curve $C$ in $X$
around $\infty$, and two points 
\[
   P^\pm\in X^\pm\cap C,\quad X^\pm=\{x\in X\mid \pm\Im x>0\},
\]
such that the image surface has cuspidal edges only along
$\mathscr{S}(C-\{P^+,P^-\})$, and has swallowtails 
only at $\mathscr{S}(P^\pm)$. 
We made our best to visualize the image surfaces; 
we often show part of the surfaces, several copies of the images of
$X^\pm$, since each of the images of the three intervals
$(-\infty,0),(0,1)$ and $(1,+\infty)$ lies on a totally geodesic surface
in $\HH^3$.

In a computational aspect of this visualization, 
we use the composition of
the hyperbolic Schwarz map $\mathscr S$ and the inverse of the Schwarz map $s$,
$\Phi=\mathscr{S}\circ s^{-1}$,
especially when the inverse of the Schwarz map is single-valued globally;
refer to Section~\ref{sec:schwarz}.
This choice is very useful, because the 
inverse map is often given explicitly as an automorphic function for the
monodromy group acting properly discontinuously on the image of the
Schwarz map. Moreover, in one of the cases where we treat the lambda
function for drawing pictures, it is indispensable,
because we have a series that converges very fast.

In the forthcoming papers we mentioned, we inroduce the derived Schwarz map,
investigate an associated parallel family of flat fronts, and
study confluence of swallowtail singularities.
Basic ingredients of the hypergeometric function and its Schwarz map can be 
found in \cite{IKSY} and \cite{Yo}.

\noindent {\sl Acknowledgment}: After the authors 
posted the first manuscript in arXiv, 
Professor D. Dumas informed them about the paper
\cite{Ep}, which introduced the notion of hyperbolic Gauss map
to nonsingular surfaces in $\HH^3$, and made some useful comments. 
The referee also gave them helpfful comments on completing the manuscript. 
The authors would like to thank both of them.

\section{Preliminaries}\label{sec:prelim}
\subsection{Models of the hyperbolic $3$-space}\label{subsec:model}
The hyperbolic $3$-space $\HH^3=\Her^+(2)/\R^+$ can be identified 
with the upper half-space $\C\times\R^+$ as
\begin{align*}
 \C\times\R^+ \ni (z,t)&\longmapsto 
        \left(\begin{array}{cc}t^2+|z|^2
	      &\bar z\\z&1\end{array}\right)\in\Her^+(2),\\
	\left(\begin{array}{cc}h&\bar w\\w&k\end{array}\right)
          \in\Her^+(2)
                         &\longmapsto
 \C\times\R^+ \ni \dfrac1k\biggl(w,\sqrt{hk-|w|^2}\biggr).
\end{align*}
It can be also identified with a subvariety 
\[
       L_1=\{x_0^2-x_1^2-x_2^2-x_3^2=1\}
\]
of the Lorentz-Minkowski $4$-space
\[
  L(+,-,-,-)=
   \bigl\{(x_0,x_1,x_2,x_3)\in\R^4\mid x_0^2-x_1^2-x_2^2-x_3^2>0,
   x_0>0\bigr\}
\]
by
\[
  \Her^+(2)\ni 
    \left(\begin{array}{cc}h&\bar w\\w&k\end{array}\right)
   \longmapsto
  \dfrac1{2\sqrt{hk-|w|^2}}
     \left(h+k,w+\bar w,\frac{w-\bar w}i,h-k\right)
                \in L_1
\]
and with the Poincar\'e ball 
\[
  B_3=\{(x_1,x_2,x_3)\in\R^3\mid x_1^2+x_2^2+x_3^2<1\},
\]
by
\[
   L_1\ni(x_0,x_1,x_2,x_3)\longmapsto\frac1{1+x_0}(x_1,x_2,x_3)\in 
   B_3.
\]
We use these models according to convenience.

\subsection{Local exponents and transformation into the $SL$-form}
The local exponents of the equation \ref{eq:hgde} at $0,1$ and $\infty$ are 
given as $\{0,1-c\},\{0,1-a-b\}$ and $\{a,b\}$, respectively.
Denote the differences of the local exponents by
\begin{equation}\label{eq:mu}
 \mu_0=1-c,\quad\mu_1=c-a-b,\quad\mu_\infty=b-a.
\end{equation}
The equation \ref{eq:hgde} transforms into the $SL$-form 
\eqref{eq:sl-form} with
\[
    q=-\frac14\left\{\frac{1-\mu_0^2}{x^2}
       +\frac{1-\mu_1^2}{(1-x)^2}
       +\frac{1+\mu_\infty^2-\mu_0^2-\mu_1^2}{x(1-x)}\right\},
\] 
by the projective change of the unknown
\[
   u\longmapsto \sqrt{x^c(1-x)^{a+b+1-c}}\ u.
\] 
Unless otherwise stated, we always take a pair $(u_0,u_1)$ of linearly 
independent solutions of \eqref{eq:sl-form} satisfying $u_0u'_1-u'_0u_1=1$,
and set
\[
  U=\left(\begin{array}{cc}u_0&u_0'\\u_1&u_1'\end{array}\right).
\]
\subsection{Monodromy group}
The group of isometries of $\HH^3$ is generated by the orientation 
preserving ones 
\[
   H\longmapsto P H \transpose{\overline P},\qquad H\in\HH^3,
      \quad P\in GL_2(\C),
\]
and the reversing one $H\to{}^tH$.

Let $\{u_0,u_1\}$ be a pair of linearly independent solutions of 
\eqref{eq:sl-form} and $\{v_0,v_1\}$ another such pair. 
Put
\[
   U=\left(\begin{array}{cc}u_0&u_0'\\u_1&u_1'\end{array}\right)
   \quad\mbox{and}\quad
   V=\left(\begin{array}{cc}v_0&v_0'\\v_1&v_1'\end{array}\right).
\]
Then there is a non-singular matrix, say $P$, such that $U=PV$ and so
that
\[
   U\transpose{\overline U}= 
   PV\transpose{\overline  V}\transpose{\overline P}.
\]
Thus the hyperbolic Schwarz map
\begin{equation}\label{eq:hsm}
   \mathscr{S}\colon{}
      X\ni x\longmapsto H(x)=
            U(x)\transpose{\overline U}(x)=
	    \left(
	    \begin{array}{cc}
	     |u_0|^2+|u_0'|^2& 
	      u_1\bar u_0+u_1'\bar u_0'\\
	    \bar u_1u_0+\bar u_1'u_0'
	    &|u_1|^2+|u_1'|^2\end{array}
	    \right)\in\HH^3
\end{equation}
is determined by the system up to orientation preserving automorphisms.
The monodromy group $\Mon(a,b,c)$ with respect to $U$ acts 
naturally on $\HH^3$ by
\[
   H\longmapsto MH\transpose{\overline M}, \quad M\in \Mon(a,b,c).
\]
Note that the hyperbolic Schwarz map to the upper half-space model is given by
\[
  X\ni x\longmapsto 
    \dfrac{\Bigl(u_0(x)\bar u_1(x)+u'_0(x)\bar u'_1(x),\ 1\Bigr)}
               {|u_1(x)|^2+|u'_1(x)|^2}\in \C\times\R^+.
\]

\subsection{Singularities of fronts}
A smooth map $f$ from a domain $U\subset \R^2$ to a Riemannian
$3$-manifold $N^3$ is called a {\em front\/} if 
there exists a unit vector field $\nu\colon{}U\to T_1N$ 
along the map $f$ such that $df$ and $\nu$ are perpendicular
and the map $\nu\colon{}U\to T_1N$ is an {\em immersion},
where $T_1N$ is the unit tangent bundle of $N$.
We call $\nu$ the {\em unit normal vector field\/} of $f$.
Not that, if we identify $T_1N$ with the unit cotangent bundle
$T_1N^*$, 
the condition $df\perp \nu$ is equivalent to  the corresponding map
$L\colon{}U\to T_1^*N$ to be Legendrian with respect to 
the canonical contact structure $T_1^*N$.
A point $x\in U$ is called a {\em singular point\/} of $f$
if the rank $df$ is less than $2$ at $x$.
It is well-known that generic singularities of fronts are 
{\em cuspidal edges\/} and {\em swallowtails} \cite{AZV}.
In this section, we roughly review these types of singularities.
  General criteria for fronts to be cuspidal edges or swallowtails
  are given in \cite{KRSUY}.

\subsubsection{$(2,3)$-cusp and cuspidal edges}
\begin{figure}
\begin{center}
\includegraphics[width=10cm]{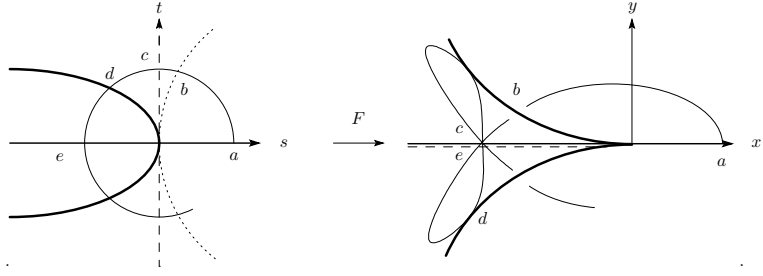}
\end{center}
\caption{Image under the map $F$}
\label{fig:F}
\end{figure}

Recall that the cubic equation
$t^3+xt-y=0$ in $t$ with real parameters $(x,y)$ has three distinct real
roots if and only if its discriminant $27y^2+4x^3$
is negative.
Consider the map
\[
   F:\R^2\ni (s,t)\longmapsto (x,y)=(s-t^2,st)\in\R^2,
\]
whose Jacobian is equal to $s+2t^2.$ 
The image of the (smooth) curve $C:s+2t^2=0$ under $F$ is a curve 
with a cusp of $(2,3)$-type, and is given by $F(C):27y^2+4x^3=0$.
Note that $F$ folds the $t$-axis to the negative half of the $x$-axis, 
and that the inverse image of $F(C)$ consists of $C$ and a curve tangent 
to $C$ at the origin, indeed we have 
\[
  27y^2+4x^3|_{x=s-t^2,\; y=st}=(s+2t^2)^2(4s-t^2).
\]
The hemicircle centered at the origin in the $(s,t)$-space is mapped 
by $F$ as is shown in Figure~\ref{fig:F}.
\par\smallskip\noindent
When a $(2,3)$ cusp traveling along a curve transversal to
$\R^2\subset\R^3$,
the locus of the singularity consists of {\em cuspidal edges}.
Precisely speaking, $p\in U$ is a {\em cuspidal edge} of a front
$f\colon{}U\to\R^3$ if there exists local diffeomorphisms
$\psi$ and $\Psi$ of $(U,p)$ and $(\R^3,f(p))$ such that
$\Psi\circ f\circ\psi(u,v)=(u^2,u^3,v)=:f_c$.
In other words, the germ of the map $f$ at $p$ is locally
{\em $A$-equivalent\/} to $f_c$.
\subsubsection{Swallowtails}
\begin{figure}
\begin{center}
\includegraphics[width=10cm]{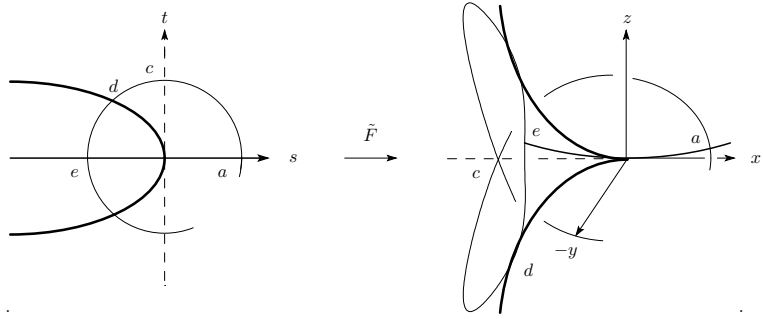}
\end{center}
\caption{Swallowtail: Image under the map $\tilde{F}$}
\label{tildeF}
\end{figure}

Consider the map
\[
  \widetilde{F}:\R^2\ni (s,t)\longmapsto
  (x,y,z)=(s-t^2,st,s^2-4st^2)\in\R^3.
\]
This map is singular (rank of the differential is not full) 
along the curve $C$, and the image of the point $(-2t^2,t)\in C$ is
given as $(-3t^2,-2t^3, 12t^4)$.
The hemicircle centered at the origin in the $(s,t)$-space 
is mapped by $\widetilde{F}$ as is shown in Figure~\ref{tildeF}. 
The image surface has three kinds of singularities
\begin{enumerate}
 \item {\it Cuspidal edges} along $\tilde F(C)-\{(0,0,0)\}$,
 \item A {\it swallowtail} at $\{(0,0,0)\}$,
 \item {\it Self-intersection} along the image of the $t$-axis.
\end{enumerate}
Here, by definition, a swallowtail is a singular point 
of a differential map $f\colon{}U\to\R^3$, which is 
$A$-equivalent to $\widetilde F(s,t)$.
Another canonical form of the swallowtail is 
$$f_s(u,v)=(3u^4+u^2v,4u^3+2uv,v),$$
which is $A$-equivalent
to $\widetilde F$ as
  $ f_s(u,v) = \Psi\circ \widetilde F\circ \psi(u,v),$ where 
$$ 
\psi(u,v)=(2v+4u^2,2u), \quad
 \Psi(x,y,z)=\left(
             \frac{-z+4x^2}{16}, \frac{y}{2},\frac{x}{2}
            \right).
$$
\section{Use of the Schwarz map}\label{sec:schwarz}
Let $u$ and $v$ be solutions of the equation \eqref{eq:sl-form}
such that $uv'-vu'=1$.
The Schwarz map is defined as 
$X\ni x\mapsto z=u(x)/v(x)\in Z$,
which is the {\em hyperbolic Gauss map\/} 
(see Section~\ref{sec:schwarz})
of the hyperbolic Schwarz map $\mathscr{S}$ as in \eqref{eq:hsm}.
It is convenient to study the hyperbolic Schwarz map 
\eqref{eq:hsm} by regarding $z$ as  variable.

Especially when the inverse of the Schwarz map is single-valued 
globally, this choice of variable is very useful, because the 
inverse map is often given explicitly as an automorphic function for the
monodromy group acting properly discontinuously on the image of the
Schwarz map.
In particular, the equation \eqref{eq:sl-form-mat} is written as
\[
   \frac{dU}{dz} 
      = U \begin{pmatrix}
	   0 & \theta \\
	  \omega & 0 \end{pmatrix}, \qquad
      \text{where}\quad
      \theta = q \frac{dx}{dz}\quad \omega=\frac{dx}{dz}.
\]
Then by the representation formula in \cite{KUY}, 
the solution $U$ is written by $\omega$, the hyperbolic Gauss map
(i.e., the Schwarz map) $z$ and their derivatives:
\begin{equation}\label{eq:kuy}
  U=i\dfrac1{\sqrt{\dot{x}}}\left(
  \begin{array}{cc}
   z\dot{x}&1+\dfrac{z}{2}\dfrac{\ddot{x}}{\dot{x}}\\[6pt]
    \dot{x}& \dfrac{1}{2}\dfrac{\ddot{x}}{\dot{x}}
  \end{array}
\right),
\end{equation}
where $\dot{~} = d/dz$.
Here, we summarize the way to show the formula:
Since $z'(:=dz/dx)=-1/v^2$ and $\ddot{x}=d^2x/dz^2$, we have 
\[
    v=i\sqrt\frac1{z'}=i\sqrt{\dot{x}},\quad u=vz,
\]
and 
\[
  v'=\dfrac{dv}{dx}=
  \dfrac{dv}{dz}\dfrac{dz}{dx}
  =\dfrac{i}2(\dot{x})^{-3/2}\ddot{x},\qquad 
  u'=i\dfrac1{\sqrt{\dot{x}}}+z\dfrac{i}2(\dot{x})^{-3/2}\ddot{x}.
\]
So, we have \eqref{eq:kuy}
and 
\begin{equation}\label{eq:H}
  H=U\transpose{\overline U}=\dfrac1{|\dot{x}|}\left(
 \begin{array}{cc}
  |z|^2|\dot{x}|^2+\left|1+\dfrac{z}{2}\,\dfrac{\ddot{x}}{\dot{x}}\right|^2
  &
  z|\dot{x}|^2
  +\dfrac{1}{2}\left(1+\dfrac{z}{2}\,\dfrac{\ddot{x}}{\dot{x}}\right)
   \dfrac{\overline{\ddot{x}}}{\overline{\dot{x}}} \\[4mm]
   \bar{z} |\dot{x}|^2+\dfrac{1}{2}
   \left(
    1+\dfrac{\bar{z}}{2}\,
    \dfrac{\overline{\ddot{x}}}{\overline{\dot{x}}}
   \right)\dfrac{\ddot{x}}{\dot{x}} 
   &|\dot{x}|^2  
   +\dfrac{1}{4}\,\left|\dfrac{\ddot{x}}{\dot{x}}\right|^2
 \end{array}
 \right).
\end{equation}
When the (projective) monodromy group of the equation
\ref{eq:hgde} is a polyhedral group or a Fuchsian triangle group, 
there is a set of real parameters $(\bar a,\bar b,\bar c)$ 
such that $\bar a-a, \bar b-b, \bar c-c\in\Z$, 
and that the Schwarz map of $E(\bar a,\bar b,\bar c)$ has 
the single-valued inverse. 
Such equations are said to be {\it standard}. 
The equation $E(a,b,c)$ is standard if $a,b,c\in\R$ satisfy 
\[
  k_0:=\frac1{|\mu_0|},\quad 
  k_1:=\frac1{|\mu_1|}, \quad 
  k_\infty:=\frac1{|\mu_\infty|}\quad\in\quad \{2,3,\dots,\infty\}.
\]
Though it is a challenging problem to study transformations of the hyperbolic 
Schwarz maps of standard equations to the general equations, we 
study only standard ones in this paper.
\section{Singularities of hyperbolic Schwarz maps}\label{sec:sing}
Since the equation \eqref{eq:sl-form} has singularities 
at $0$, $1$ and $\infty$, the corresponding hyperbolic Schwarz map 
$\mathscr{S}$ has singularities at these points.
In terms of flat fronts in $\HH^3$, they are considered as {\em ends\/}
of the surface.
On the other hand, the map $\mathscr{S}$
may not be an immersion at  $x\in X$, even
if  $x$ is not a singular point of \eqref{eq:sl-form}.
In other words, $x$ is a singular point of the front 
$\mathscr{S}\colon{}X\to \HH^3$.

In this section, we analyze properties of these singular points of
the hyperbolic Schwarz maps.
\subsection{Singularities on $X$}
As we have seen in the introduction, the hyperbolic Schwarz map
$\mathscr{S}\colon{}X=\C-\{0,1\}\to \HH^3$ 
can be considered as a flat front in the sense of \cite{KUY,KRSUY}.
Thus, as a corollary of Theorem~1.1 in \cite{KRSUY}, we have 
\begin{lemma}\label{lem:criteria}
 \begin{enumerate}
  \item  A point $p\in X$ is a singular point of the hyperbolic Schwarz
	 map $\mathscr{S}$ if and only if $|q(p)|=1$,
  \item  a singular point $x\in X$ of $H$ is A-equivalent to 
	 the cuspidal edge if and only if 
	 \[
	    q'(x) \neq 0 \qquad \text{and}\qquad
	    q^3(x)\bar q'(x) - q'(x)\neq 0,
	 \]
  \item  and singular point $x\in X$ of $H$ is A-equivalent to 
	 the swallowtail if and only if 
	 \begin{align*}
	    q'(x) \neq 0 &\qquad 
	    q^3(x)\bar q'(x) - q'(x)= 0,\\\text{and}\quad
	   & \Re
	   \left\{\frac{1}{q}\left(
	        \left(\frac{q'(x)}{q(x)}\right)' -
	        \frac{1}{2}\left(
	                    \frac{q'(x)}{q(x)}
                           \right)^2
	                     \right)\right\}\neq 0.
	 \end{align*}
 \end{enumerate}
\end{lemma}
We apply Lemma~\ref{lem:criteria} to the hypergeometric equation.
Using $\mu_0$, $\mu_1$ and $\mu_{\infty}$ as in \eqref{eq:mu},
the coefficient of the hypergeometric equation \eqref{eq:sl-form} is written as
\begin{equation}\label{eq:q}
  q=-\frac14\left(\frac{1-\mu_0^2}{x^2}+\frac{1-\mu_1^2}{(1-x)^2}
 +\frac{1+\mu_\infty^2-\mu_0^2-\mu_1^2}{x(1-x)}\right)
 =:\frac{-Q}{4x^2(1-x)^2},
\end{equation}
where
\begin{equation}\label{eq:Q}
  Q=1-\mu_0^2+(\mu_\infty^2+\mu_0^2-\mu_1^2-1)x+(1-\mu_\infty^2)x^2.
\end{equation}
Hence $x\in X$ is a singular point if and only if
\begin{equation}\label{eq:singQ}
   |Q| = 4 |x^2(1-x)^2|.
\end{equation}
Define $R$ by
\begin{equation}\label{eq:R}
  q'=-\frac{Q'x(1-x)-2Q(1-2x)}{4x^3(1-x)^3}=:\frac{-R}{4x^3(1-x)^3}.
\end{equation}
Then  we have 
\[
   q^3(x)\bar q'(x)-q'(x)=\frac{Q^3}{4^3x^6(1-x)^6}\cdot
  \frac{\overline R}{4\bar x^3(1-\bar x)^3}+\frac{ R}{4x^3(1- x)^3}.
\]
So, 
the condition $q^3(x)\bar q'(x)-q'(x)=0$ is equivalent to 
`$Q_{}^3 \overline{R}{}^2$ is real non-positive' 
under the condition \eqref{eq:singQ}. 
Therefore, a singular point $x$ is a cuspidal edge if and only if 
\eqref{eq:singQ} is satisfied and
\begin{equation}\label{eq:cuspQ}
   Q_{}^3 \overline{R}{}^2\quad \mbox{is not a non-positive real number.}
\end{equation}
Moreover, a singular point $x$ is a swallowtail if and only if
\begin{equation}\label{eq:swtlQ}
    Q_{}^3 \overline{R}{}^2 \mbox{ is real non-positive, and }
        \Re\left(
	  2|R|^4-x(1-x)(2R'Q-RQ')\overline R^2
        \right)\neq 0
\end{equation}
hold, where $'=d/dx$.
In fact, since $(q'/q)=R/\bigl(x(1-x)Q\bigr)$, we have
\begin{align*}
 \frac{1}{q}&\left(
            \left(\frac{q'}{q}\right)'-\frac{1}{2}\left(\frac{q'}{q}\right)^2
            \right)\\
  &= -\frac{2}{Q^3}\left(
       2(1-2x)RQ+2x(1-x)(R'Q-RQ')-R^2
     \right)\\
  &=-\frac{2}{Q^3}(-2R^2+x(1-x)(2R'Q-RQ'))\\
  &= \frac{2R_{}^2\overline{Q}^3}{|R_{}^2\overline Q^3|^2}
         \left(2|R|^4-x(1-x)(2R'Q-RQ')\overline R^2\right).
\end{align*}

\subsection{At a singular point of the equation \ref{eq:hgde}}
 In this subsection, we assume the parameters $a$, $b$ and $c$ are
 real.
 Since $q$ has poles of order $2$ at $0$, $1$ and $\infty$,
 $|q|\neq 1$ in a neighborhood of the singularities of the equation
 \ref{eq:hgde}.
 Here, we study the behavior of $X$ around these points.
 If $X$ were single-valued on a neighborhood of the end, 
 the following calculations are essentially similar 
 to those in \cite{GMM}, in which asymptotic behavior 
 of the end of flat fronts is investigated.

 Around, for example, $x=0$, the Schwarz map has the expression 
 $z=x^{|1-c|}(1+O(x))$. 
 So we may assume that the inverse map has 
 the expression $x=z^\alpha\bigl(1+O(z)\bigr)$ for some {\it real\/}
 constant $\alpha(>0)$.
 Since 
\[
   \dot{x}=\alpha z^{\alpha-1}(1+O(z)),\quad 
       \frac{\ddot{x}}{\dot{x}}=\frac{\alpha-1}z(1+O(z)),
\]
the principal part of the matrix $U$ is given by
\[
   P:=\frac{i}{\sqrt{\alpha z^{\alpha+1}}}
   \left(
   \begin{array}{cc}
    \alpha z^{\alpha+1}
     &\left(1+\dfrac{\alpha-1}2\right)z\\[3mm]
     \alpha z^\alpha&\dfrac{\alpha-1}2
   \end{array}
   \right).
\]
We have
\[
    P \transpose{\overline P}:=\frac1{|\alpha z^{\alpha+1}|}
    \left(
    \begin{array}{cc} 
     \qquad **\qquad 
      &\left(\dfrac{\alpha^2-1}4\right)z+|\alpha z^\alpha|^2z\\[3mm]
      **&\left(\dfrac{\alpha-1}2\right)^2+|\alpha z^\alpha|^2
    \end{array}
    \right).
\]
Thus the hyperbolic Schwarz map $\mathscr{S}$ extends to the point $z=0$ and 
to the boundary of $\HH^3$. Its image is nonsingular at $\mathscr{S}(0)$, 
and is tangent to the boundary at this point. 
\section{Hyperbolic Schwarz maps}
When the monodromy group of the equation $E(a,b,c)$ is a finite group or a
typical Fuchsian group, we study the singularities of the hyperbolic 
Schwarz map, and visualize the image surface.
\subsection{Finite (polyhedral) monodromy groups}
We first recall fundamental facts about the polyhedral groups and their
invariants basically following \cite{Kl}.
\subsubsection{Basic data}
Let the triple $(k_0,k_1,k_\infty)$ be one of 
\[
   (2,2,n)\ (n=1,2,\dots),\quad(2,3,3),\quad(2,3,4),\quad(2,3,5),
\]
in which case, the projective monodromy group is of finite order $N$: 
\[
   N=\quad2n,\quad12,\quad24,\quad60,
\]
respectively. 
Note that
\[
   \frac2N=\frac1{k_0}+\frac1{k_1}+\frac1{k_\infty}-1.
\]
For each case, we give a triplet $\{R_1,R_2,R_3\}$ of reflections 
whose mirrors bound a Schwarz triangle. 
These are tabulated in Table~\ref{tab:2}, where 
\[
   R(c,r): z\longmapsto \dfrac{c\bar z+r^2-|c|^2}{\bar z-\bar c}
\]
is the reflection with respect to the circle of radius $r>0$ centered at
$c$. The monodromy group Mon (a polyhedral group) is the group 
of even words of these three reflections.
\begin{table}
\footnotesize
\fbox{
 \begin{minipage}{12cm}
  \begin{description}
   \item[Dihedral]
      $R_1: z\mapsto \bar z$,
      $R_2: z\mapsto e^{2\pi i/n}\bar z$,
      $R_3: z\mapsto \dfrac1{\bar z}$.
   \item[Tetrahedral]\
      $R_1: z\mapsto \bar z$,
      $R_2: z\mapsto -\bar z$,
      $R_3=R(-\dfrac{1+i}{\sqrt{2}},\sqrt{2})$.
   \item[Octahedral]\
      $R_1: z\mapsto \bar z$,
      $R_2: z\mapsto i\bar z$, 
      $R_3=R(-1,\sqrt{2})$.
   \item[Icosahedral]
      $R_1: z\mapsto \bar z$,
      $R_2: z\mapsto \epsilon^2\bar z$,
   \[
     R_3=R\left(2\cos\dfrac{\pi}5,\sqrt{1+4\cos^2\dfrac{\pi}5}\right)
	      =\dfrac{-(\epsilon-\epsilon^4)\bar
	      z+(\epsilon^2-\epsilon^3)}{(\epsilon^2-\epsilon^3)\bar
	      z+(\epsilon-\epsilon^4)},\quad \epsilon=e^{2\pi i/5}.
   \]
  \end{description}
  \end{minipage}}
\caption{}\label{tab:2}
\end{table}

\noindent
The (single-valued) inverse map
\[
   s^{-1}:Z\ni z\longmapsto x\in \bar X\cong\PP^1,
\]
invariant under the action of Mon, is given as follows. 
Let $f_0(z), f_1(z)$ and $f_\infty(z)$ be the monic polynomials in $z$ 
with simple zeros exactly at the images $s(0),s(1)$ and $s(\infty)$,
respectively. If $\infty\in Z$ is not in these images, then the degrees
of these polynomials are $N/k_0,N/k_1$ and $N/k_\infty$, respectively; 
if for instance
$\infty\in s(0)$, then the degree of $f_0$ is $N/k_0-1$.
Now the inverse map $s^{-1}$ is given by
\[
   x=A_0\ \frac{f_0(z)^{k_0}}{f_\infty(z)^{k_\infty}},
\]
where $A_0$ is a constant; we also have
\[
  1-x=A_1\ \frac{f_1(z)^{k_1}}{f_\infty(z)^{k_\infty}},\quad 
  \frac{dx}{dz}=A
      \frac{f_0(z)^{k_0-1}f_1(z)^{k_1-1}}{f_\infty(z)^{k_\infty+1}},
\]
for some constants $A_1$ and $A$.
See Table~\ref{tab:1}.
\begin{table}\footnotesize
\fbox{
 \begin{minipage}{12cm}
  \begin{description}
   \item[Dihedral]
	      $(k_0,k_1,k_\infty)=(2,2,n),$ $N=2n.$
	\begin{gather*}
          A_0=\frac14,\quad A_1=-\frac14,\quad A=\frac{n}4,\\
	 f_0=z^n+1,\quad
	 f_1=z^n-1,\quad
	 f_\infty=z.
	\end{gather*}
      $f_\infty$ is of degree $1=2n/n-1$, since $\infty\in s(\infty)$, 
	      that is, $x(\infty)=\infty.$
  \item[Tetrahedral] $(k_0,k_1,k_\infty)=(2,3,3),$ $N=12$.
       \begin{gather*}	      
         A_0=-12\sqrt3,\quad A_1=1,\quad A=24\sqrt3,\\
	\begin{aligned}
	 f_0&=z(z^4+1),\\
	 f_1&=z^4+2\sqrt3z^2-1=(z^2-2+\sqrt3)(z^2+2+\sqrt3),\\
	 f_\infty&=z^4-2\sqrt3z^2-1=(z^2-2-\sqrt3)(z^2+2-\sqrt3).
	\end{aligned}
       \end{gather*}
      $f_0$ is of degree $5=12/2-1$, since $\infty\in s(0)$, that is,
      $x(\infty)=0$.
  \item[Octahedral] $(k_0,k_1,k_\infty)=(3,2,4),$ $N=24$.
      \begin{gather*}
        A_0=\frac1{108},\quad A_1=\frac{-1}{108},\quad A=\frac1{27},\\
       \begin{aligned}
	f_0&=z^8+14z^4+1=(z^4+2z^3+2z^2-2z+1)(z^4-2z^3+2z^2+2z+1),\\
	f_1&=z^{12}-33z^8-33z^4+1=(z^4+1)(z^2+2z-1)(z^2-2z-1)(z^4+6z^2+1),\\
        f_\infty&=z(z^4-1)=z(z^2+1)(z^2-1).
       \end{aligned}
       \end{gather*}
       $f_\infty$ is of degree $5=24/4-1$, since $\infty\in s(\infty)$, 
       that is, $x(\infty)=\infty$.
  \item[Icosahedral] $(k_0,k_1,k_\infty)=(3,2,5),$ $N=60$.
       \begin{gather*}
	 A_0=\frac{-1}{1728},\quad A_1=\frac1{1728},\quad
	 A=\frac{-5}{1728},\\
	\begin{aligned}
	 f_0&=z^{20}-228z^{15}+494z^{10}+228z^5+1\\
	 &=(z^4-3z^3-z^2+3z+1)(z^8-z^7+7z^6+7z^5-7z^3+7z^2+z+1)\\
	 &\hphantom{=}\times(z^8+4z^7+7z^6+2
	 z^5+15z^4-2z^3+7z^2-4z+1),\\
	 f_1&=z^{30}+522z^{25}-10005z^{20}-10005z^{10}-522z^5+1\\
	 &=(z^2+1)(z^8-z^6+z^4-z^2+1)(z^4+2z^3-6z^2-2z+1)\\
	 &\hphantom{=}\times(z^8+4z^7+17z^6+22z^5+5z^4-22z^3+17z^2-4z+1)\\
	 &\hphantom{=}\times(z^8-6z^7+17z^6-18z^5+25z^4+18z^3+17z^2+6z+1),\\
	 f_\infty&=z(z^{10}+11z^5-1)\\
	 &=z(z^2+z-1)(z^4+2z^3+4z^2+3z+1)(z^4-3z^3+4z^2-2z+1).
	\end{aligned}
	\end{gather*}
       $f_\infty$ is of degree $11=60/5-1$, since $\infty\in s(\infty)$, 
       that is, $x(\infty)=\infty.$
  \end{description}
 \end{minipage}}
\caption{}\label{tab:1}
\end{table}
\subsubsection{Dihedral cases}
\begin{figure}
\begin{center}
\begin{tabular}{cc}
 \raisebox{5.5mm}{
 \includegraphics[width=5cm]{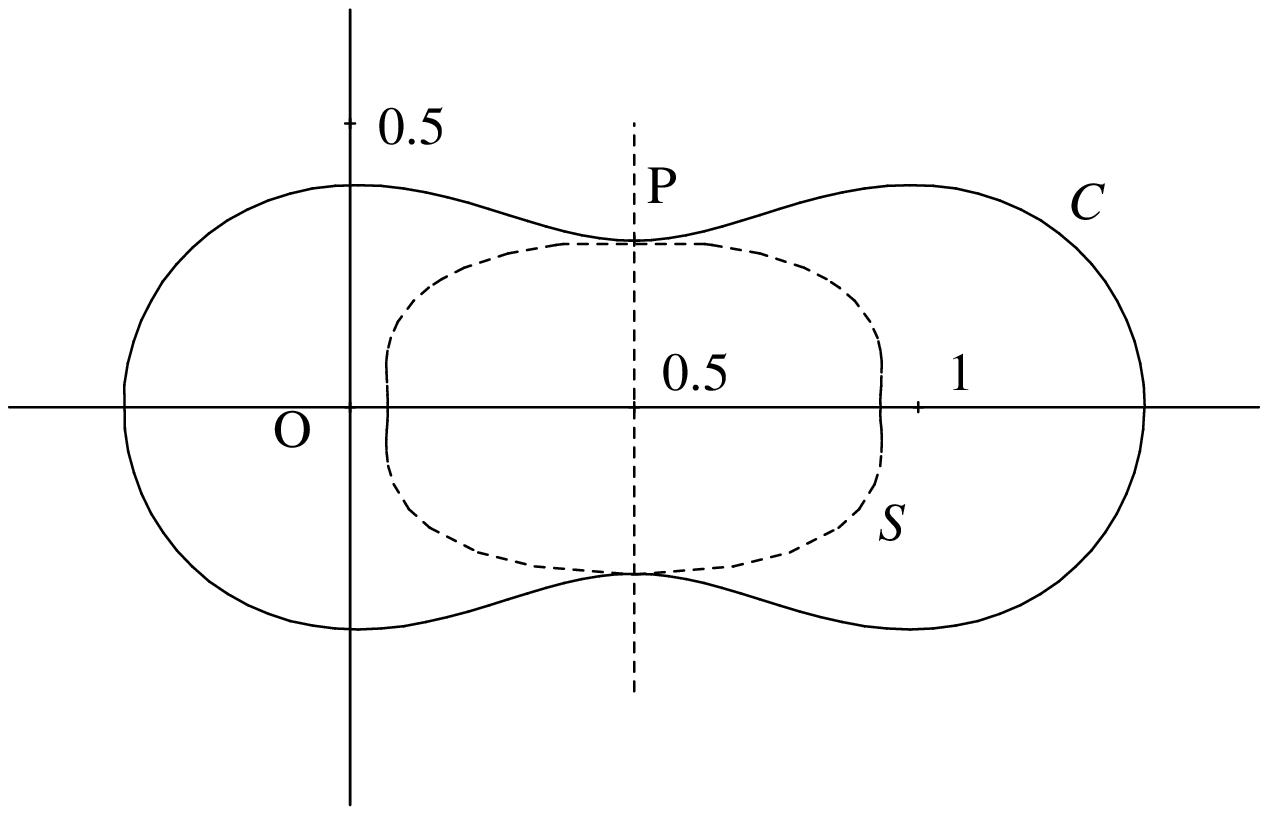}} &
 \includegraphics[width=5cm]{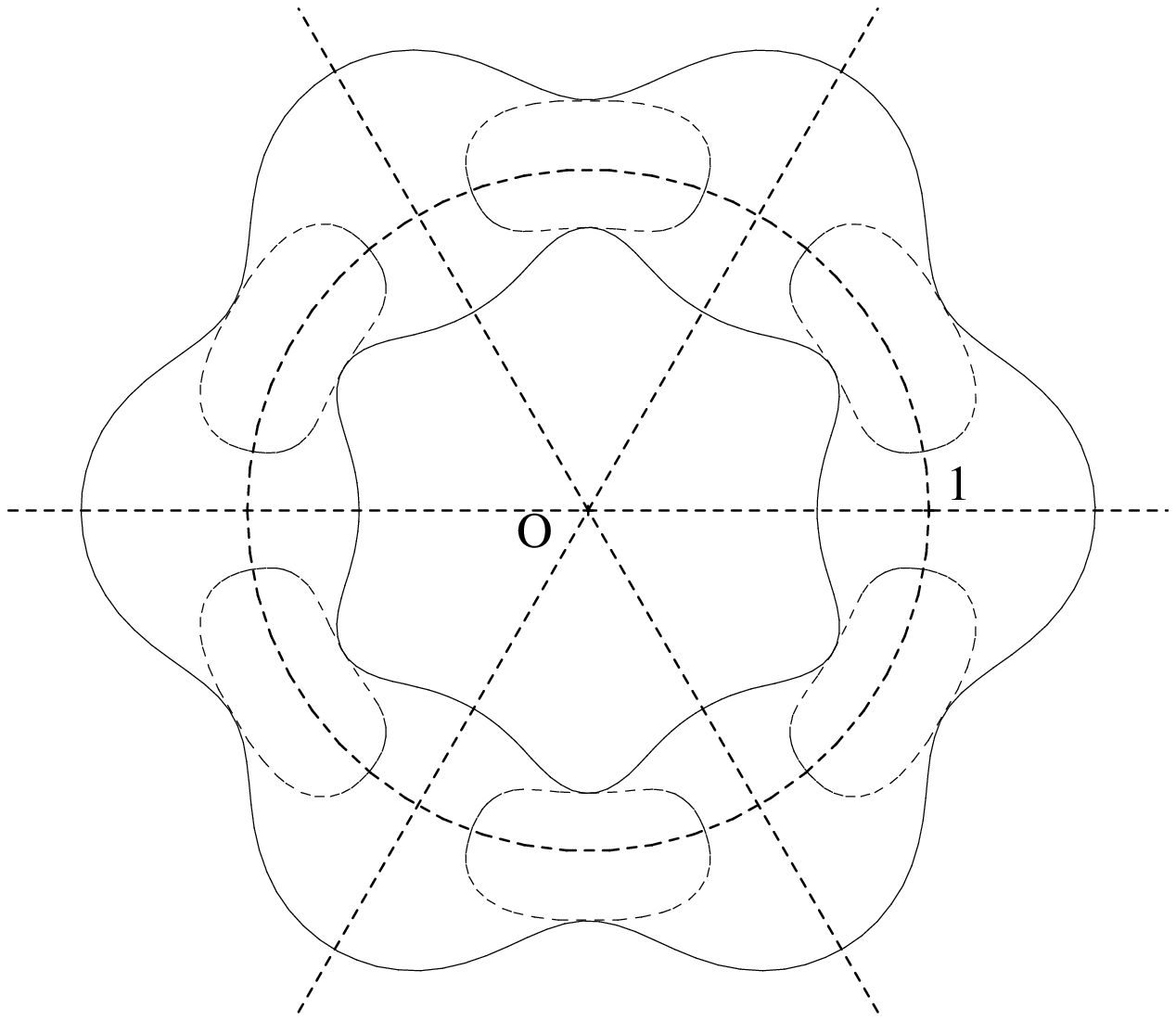} \\
 in the $x$-plane & in the $z$-plane
\end{tabular}
\end{center}
\caption{The curve $C:|Q|=4|x(1-x)|^2$, 
  when $(k_0,k_1,k_\infty)=(2,2,3)$}
\label{dihsingx}
\end{figure}

We consider a dihedral case: $(k_0,k_1,k_\infty)=(2,2,n)$, $n=3$. 
The curve $C$ in the $x$-plane defined by \eqref{eq:singQ}: $|Q|=4|x(1-x)|^2$ is 
symmetric with respect to the line $\Re(x)=1/2$ and 
has a shape of a cocoon (see Figure \ref{dihsingx} (left)).
We next study the condition \eqref{eq:swtlQ}.
The curve $\Im(Q^3\overline{R}{}^2)=0$ consists of the line $\Re x=1/2$, 
the real axis, and a curve of degree 8. 
We can prove that, on the upper half $x$-plane, 
there is a unique point satisfying the conditions \eqref{eq:singQ} and 
\eqref{eq:swtlQ},
(this point is the intersection $P$ of the curve $C$ and 
the line $\Re x=1/2$,) and that the image surface has a swallowtail
at this point, and has cuspidal edges along $\mathscr{S}(C)$ outside
$\mathscr{S}(P)$. 
We omit the proof since the computation is analogous to the case 
$(k_0,k_1,k_\infty)=(\infty,\infty, \infty)$; see
Section~\ref{subsub:fuchs}.

The curve which gives the self-intersection is tangent to $C$ at $P$,  
and crosses the real axis perpendicularly; 
this is the dotted curve in Figure \ref{dihsingx} (left), and is made as
follows. 
Since the curve is symmetric with respect
to the line $\Re x=1/2$, on each level line $\Im x=t$, we take two
points $x_1$ and $x_2$ ($\Re(x_1+x_2)=1$), compute the distance between
their images $\mathscr{S}(x_1)$ and $\mathscr{S}(x_2)$, 
and find the points that the two image points coincide. 

We substitute the inverse of the Schwarz map (cf. Table~\ref{tab:1}):
\[
   x=\frac14\ \frac{(z^n+1)^2}{z^n}, \qquad n=3
\]
into the expression \eqref{eq:H} of the hyperbolic Schwarz map,
and visualize the image surface in the Poincar\'e ball model 
explained in Section~\ref{subsec:model}.
The upper half $x$-space corresponds to a fan in the $z$-plane 
bounded by the lines with argument $0,\pi /3,2\pi /3$, 
and the unit circle (see Figure~\ref{dihsingx}(right).). 
The image $s(C)$ consists of two curves; 
the dotted curves in the figures form the pre-image
of the self-intersection.

\begin{figure}
\begin{center}
\includegraphics[width=5.5cm]{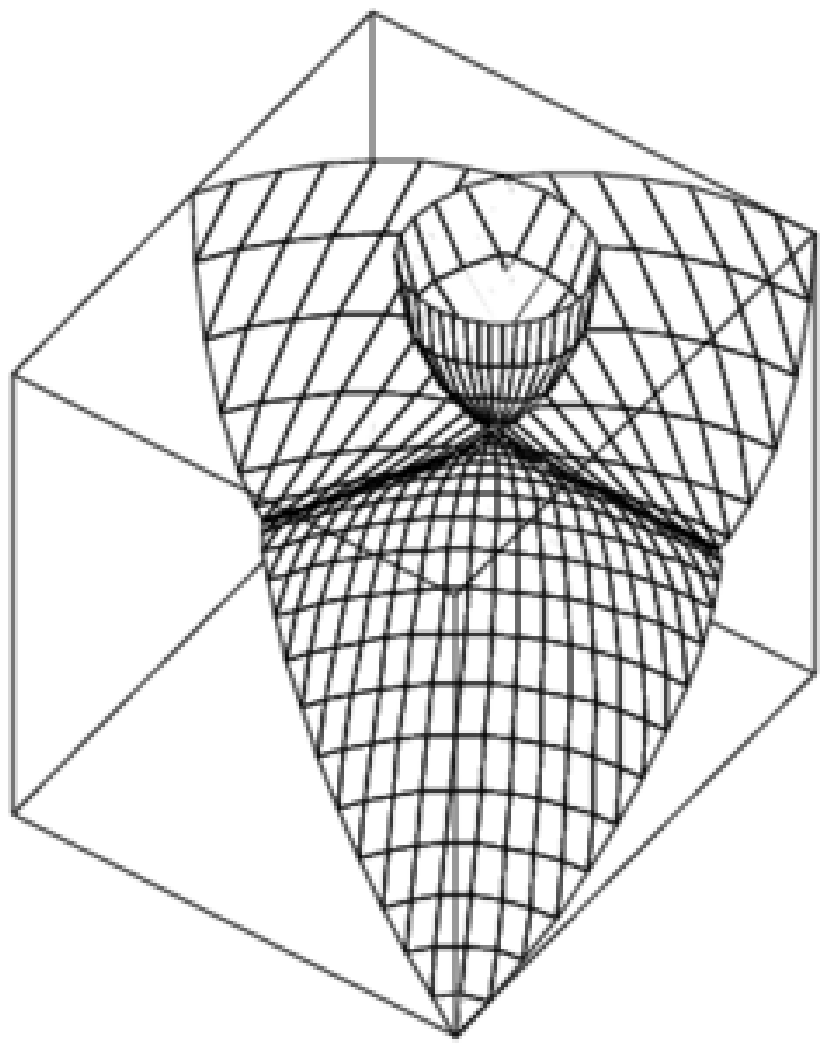}\hspace{1cm}
\includegraphics[width=5.5cm]{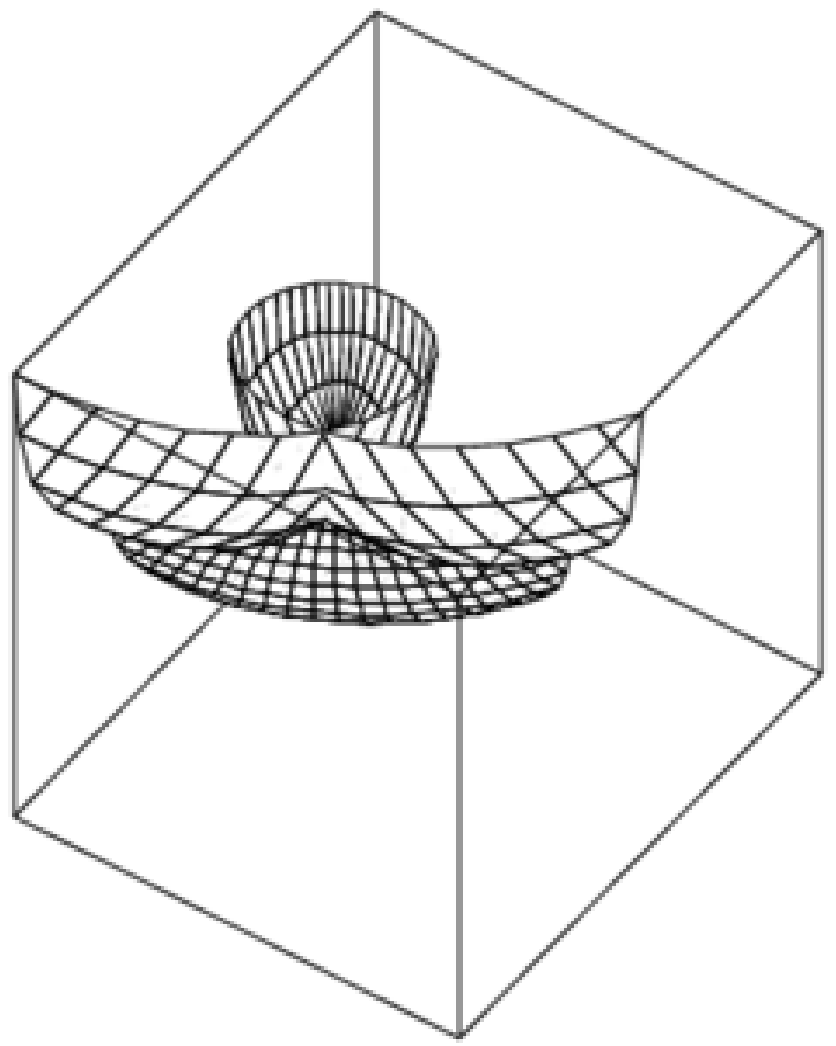}\\
   Image of a fan under $\Phi$\\
\includegraphics[width=7cm]{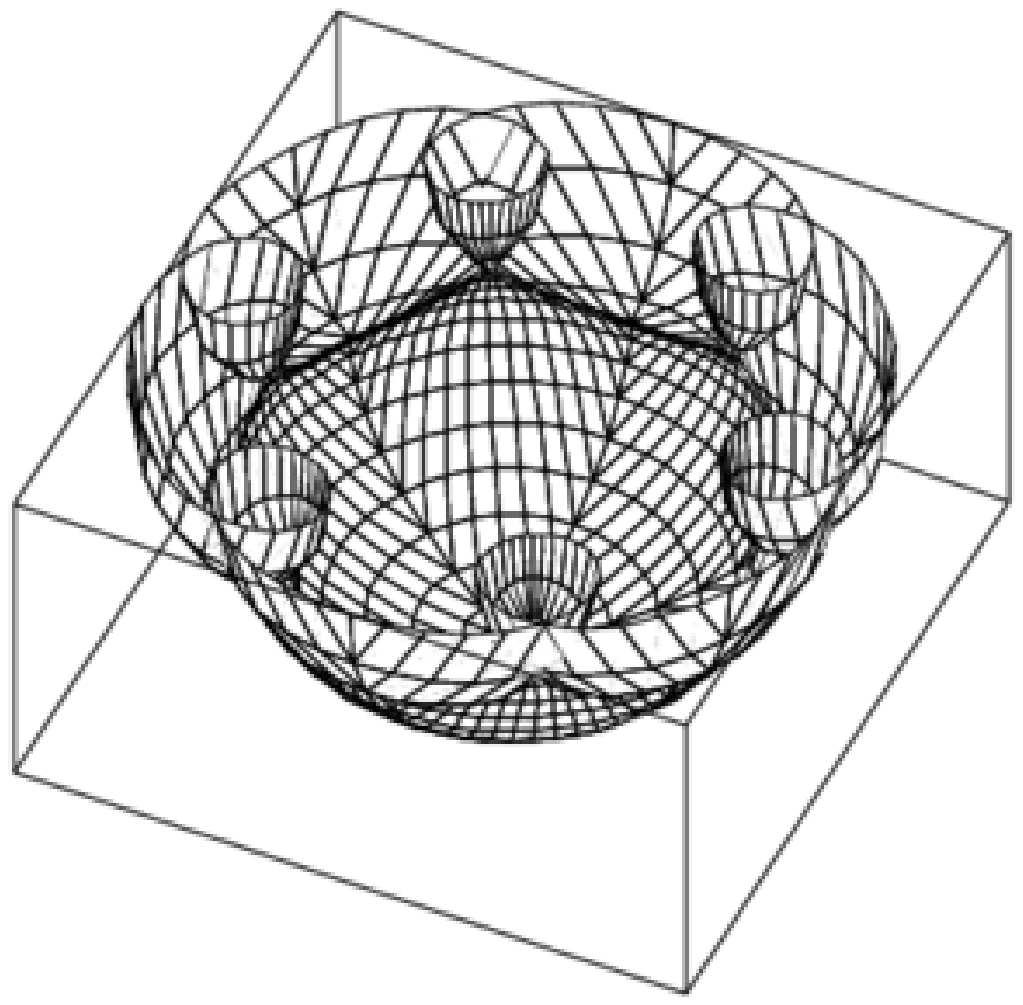}\\
Image of six fans
\end{center}
\caption{Dihedral case}
\label{dihp}
\end{figure}

Let $\Phi$ denote the hyperbolic Schwarz map in $z$-variable:
$$\Phi:=\mathscr{S}\circ s^{-1}:Z\ni z\longmapsto H(z)\in \HH^3.$$
We visualize the image of the hyperbolic Schwarz map when $n=3$.
Figure \ref{dihp}(upper left) is a view of the image 
of one fan in the $z$-plane
under $\Phi$ (equivalently, the image of upper/lower half $x$-plane under 
$\mathscr{S}$).
The cuspidal edge traverses the figure from left to right
and one swallowtail is visible in the center. The
upper right figure is the antipode of the left.
Figure \ref{dihp}(below) is a view of the image of six fans 
dividing the unit $z$-disk. To draw the images of fans with the same
accuracy, we make use of the invariance of the function $x(z)$ under the
monodromy groups.
\subsubsection{Other polyhedral cases}
\begin{figure}
\begin{center}
\begin{tabular}{cc}
 \includegraphics[height=5.5cm]{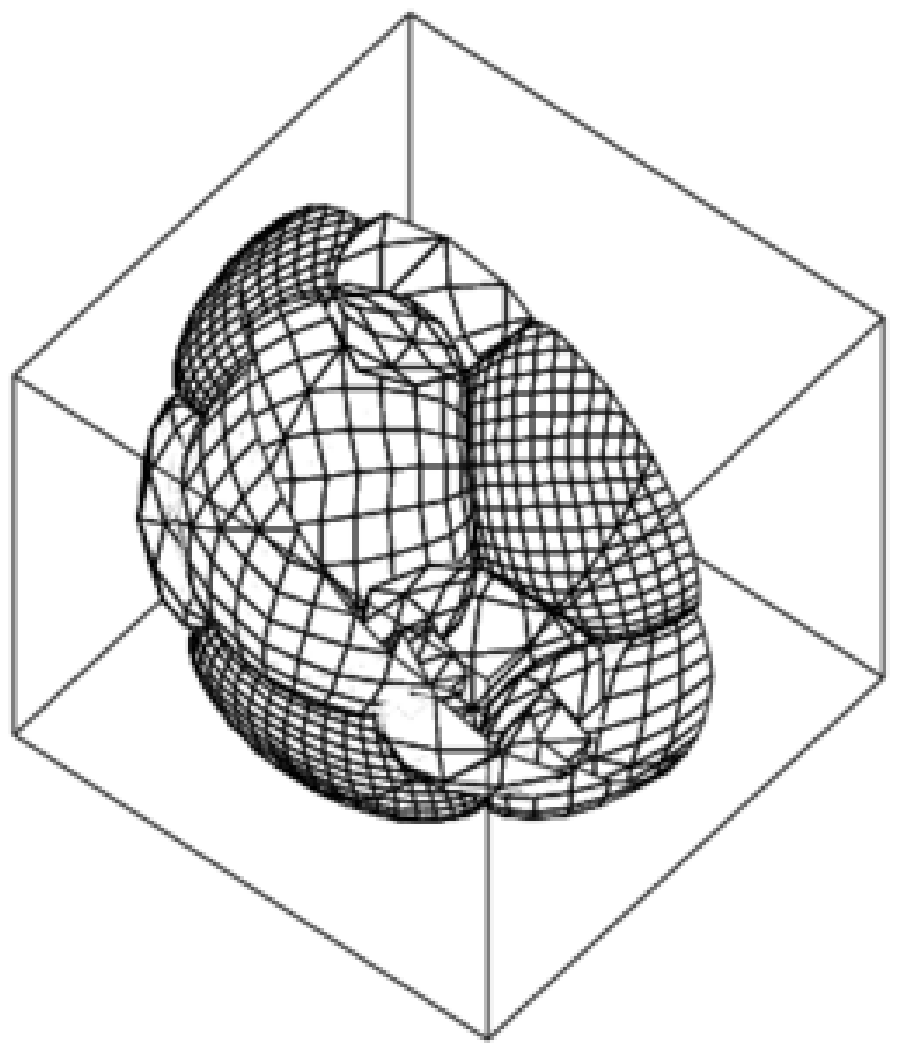} &
 \includegraphics[height=5.5cm]{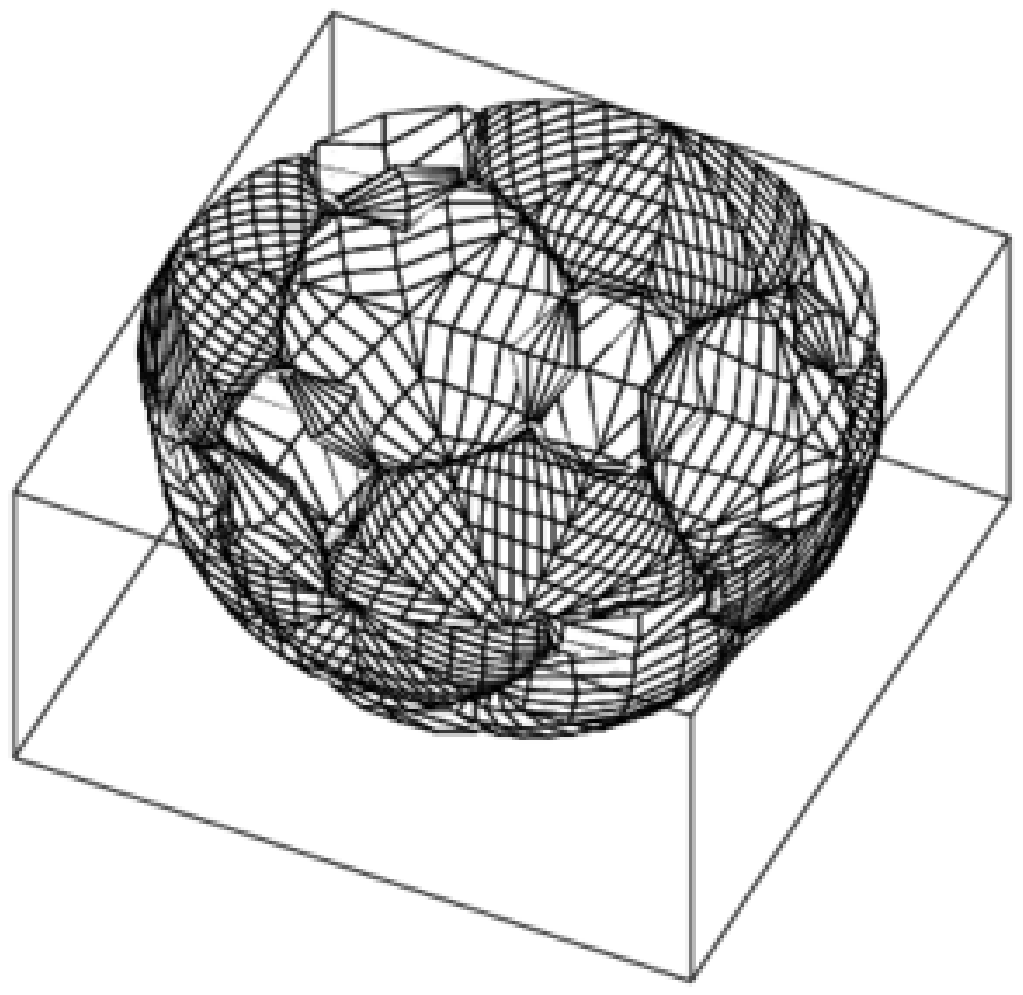} \\
 Tetrahedral&
 Octahedral\\
 Image of twelve triangles & Image of 32 triangles \\
 \multicolumn{2}{c}{
 \includegraphics[height=5.5cm]{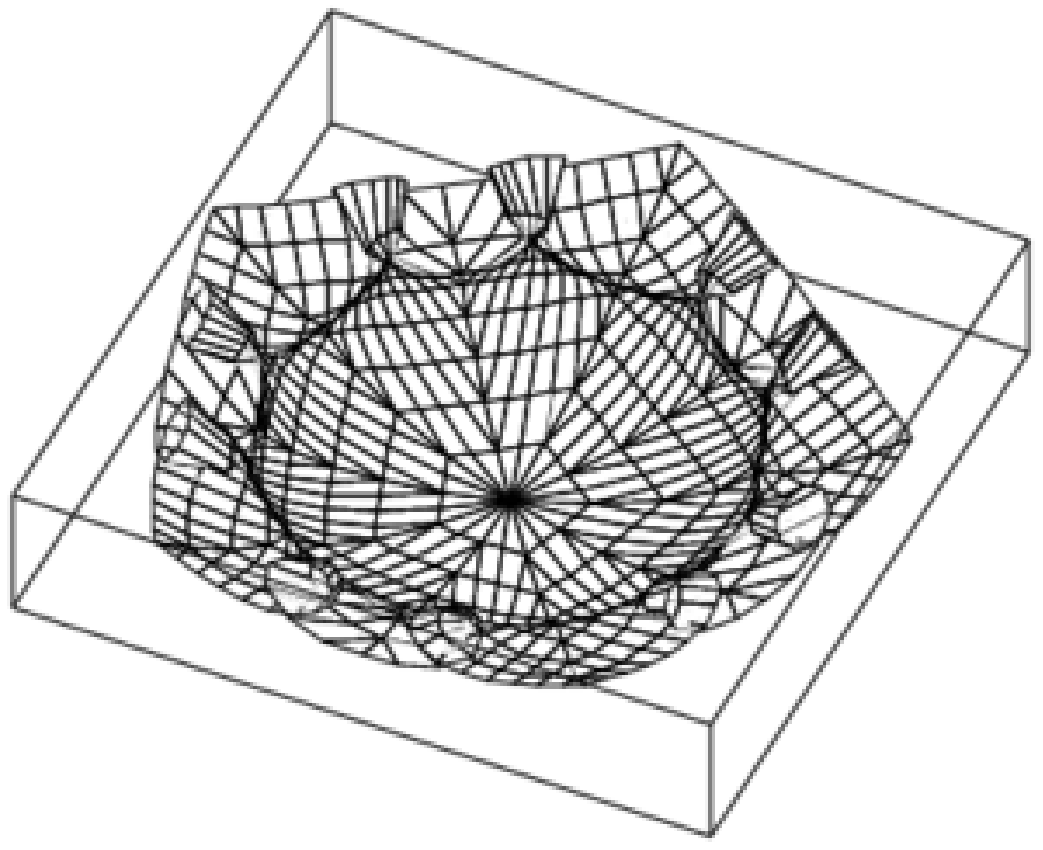} }\\
 \multicolumn{2}{c}{
 Icosahedral; Image of ten triangles }\\[5mm]
\includegraphics[height=5.5cm]{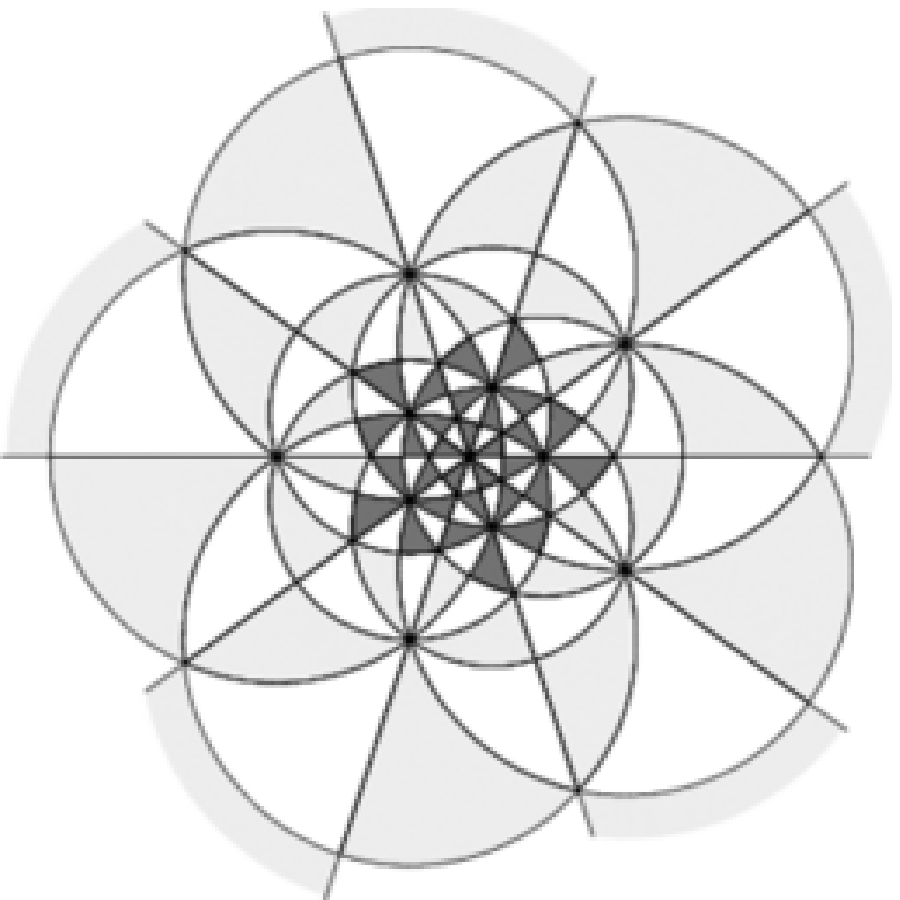} &
 \includegraphics[height=5.5cm]{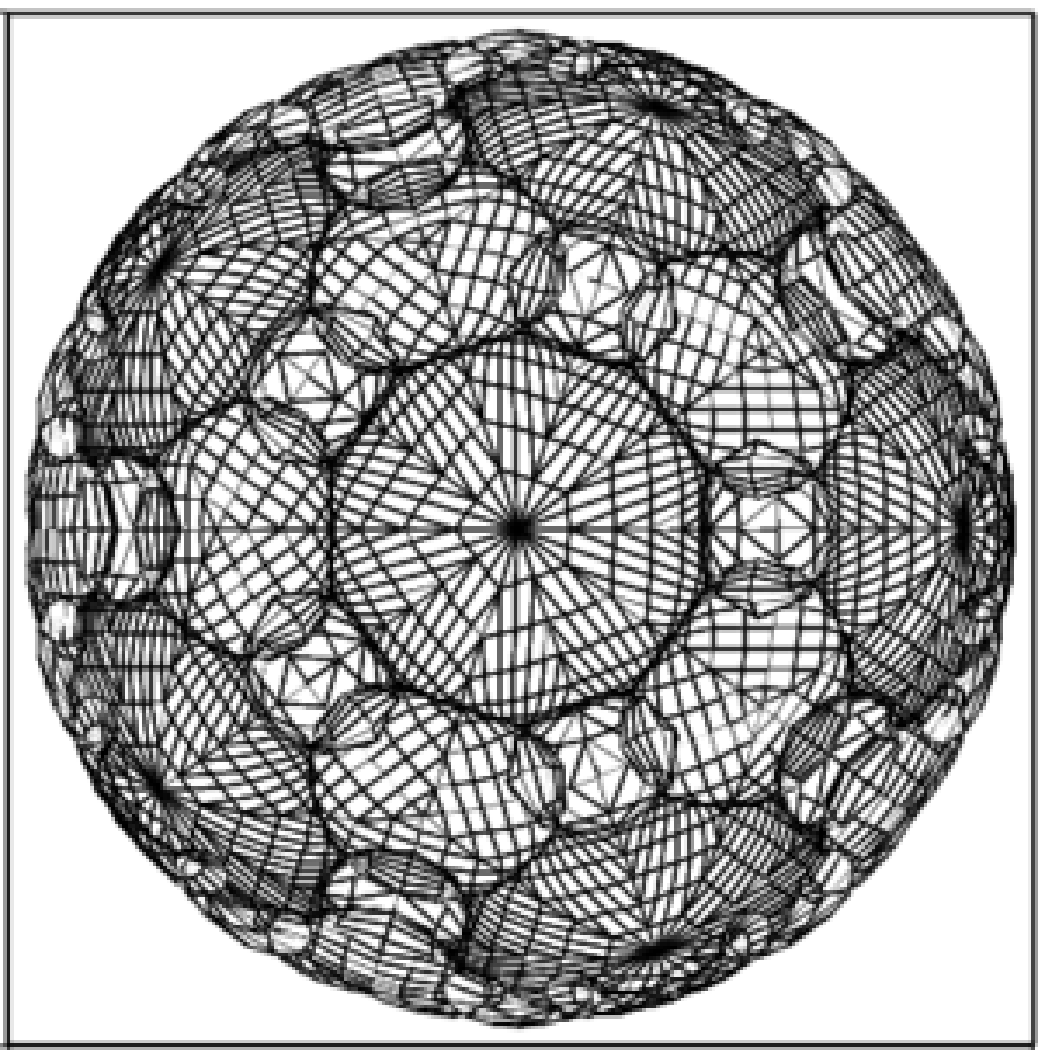} \\
 120 icosahedral triangles &
 Icosahedral\\
 and 60 central ones & Image of 60 triangles
\end{tabular}
\end{center}
\caption{Other polyhedral cases}
\label{polyhedral}
\end{figure}

For other polyhedral cases, situation is similar. 
The sphere $Z$ is divided into $2N$ triangles.
In Figure~\ref{polyhedral}, for the tetrahedral and the octahedral cases,  
the images under $\Phi$ of $N$ triangles are shown, 
and for the icosahedral case,
$2N=120$ triangles dividing the $z$-plane, the images of central ten triangles,
and that of $N=60$ triangles are shown.

\subsection{A Fuchsian monodromy group}
We study only the case $(k_0,k_1,k_\infty)=(\infty,\infty,\infty)$.
\subsubsection{Singular locus}\label{subsub:fuchs}
We find the singular locus of the image when $\mu_0=\mu_1=\mu_\infty=0$. 
We have
\[
   Q=1-x+{x}^{2},\quad R=\left (-1+2\,x\right )\left ({x}^{2}-x+2\right).
\]
The singularities lie on the image of the curve 
\[
   C: f:=16|x(1-x)|^4-|Q|^2=0.
\]
Note that this curve is symmetric with respect to the line $\Re x=1/2$.
\begin{figure}
\begin{center}
 \includegraphics[width=5.5cm]{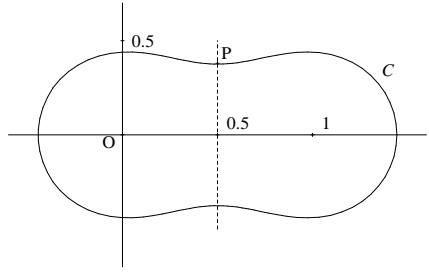}
\end{center}
\caption{The curve $C:|Q|=4|x(1-x)|^2$, 
        when $(k_0,k_1,k_\infty)=(\infty.\infty,\infty)$}
\label{lamsing}
\end{figure}

Recall that the condition \eqref{eq:swtlQ} is stated as
\[
   h:=\Im(Q^3\overline{R}{}^2)=0,\qquad \Re(Q^3\overline{R}{}^2)>0.
\]
The curve $h=0$ consists of the line $\Re x=1/2$, the real axis, 
and a curve of degree 8. We can prove that, on the upper half $x$-plane, 
there is a unique point in the intersection points of the curves $C$ 
and  $h=0$ satisfying the conditions \eqref{eq:singQ} and
\eqref{eq:swtlQ} -- 
 this point is the intersection $P$ of the curve $C$ and 
the line $\Re(x)=1/2$ --  
and that the image surface has a swallowtail singularity 
at $P$, and has cuspidal edges along $\mathscr{S}(C)$ outside $\mathscr{S}(P)$.

Actual computation proceeds as follows. The image curve has singularities 
at the image of the intersection of the curves $C$ and 
$\{h=0,\ \Re(Q^3\overline{R}{}^2)>0\}$.
We can show that there is only one such point: the intersection of $C$
and the line $\Re x=1/2.$ 

When 
$\mu_0=\mu_1=\mu_\infty=0$, the coefficient $q$ is expressed as
\[
   q=-\frac14\frac{Q}{x^2(1-x)^2},\qquad 
   Q(x)=x^2-x+1.
\]
If we put $x=s+it$, then $f:=|Q|^2-4^2|x^2(1-x)^2|^2$ is
a polynomial in $s$ and $t$ of order 8. 
If we put 
\[
   s=\frac{1}{2}+u, \quad u^2=U,\quad t^2=T,
\]
then $f$ turns out to be a polynomial $F$ in $U$ and $T$ of order $4$:
\begin{align*}
  F&= \dfrac12+\dfrac52(U-T)-5(U^2+T^2)+6TU
   +16(TU^2-T^2U)+16(U^3-T^3)\\
   &  \quad -16(U^4+T^4)-64(T^3U-TU^3)-96T^2U^2.
\end{align*}

The polynomial $R$ is expressed as
\[
   R=-4x^3(1-x)^3q'(x)=(2x-1)(x^2-x+2).
\]
The imaginary part of $Q^3\overline R^2$ has the form $t(2s-1)G$, 
where $G$ is a polynomial in $U$ and $T$ of order $4$:
\begin{align*}
  G&= \dfrac{1323}{256}+\dfrac{189}{16}(U-T)
         +\dfrac98(U^2+T^2)-\dfrac{99}{4}TU
   +11(T^3-U^3)+11(T^2U-TU^2)\\
  & \quad -5(U^4-T^4)-20(T^3U+TU^3)-30T^2U^2.
\end{align*}
Set 
\[
   G_1:=5F-16G,\quad F_1:=256F-16G_1 \quad\mbox{and}
   \quad U-T=:S,\quad UT=:V.
\]
Then we have
\[
   G_1=-\dfrac{1283}{16}+256S^3-43S^2+1024VS-\dfrac{353}{2}S+340V,
\]
which is {\it linear\/} in $V$. 
Solving $V$ from the equality $G_1=0$, 
and substituting it into $F_1=0$, we get a rational function in $S$, 
whose numerator is a polynomial in $S$ of degree 3. 
The roots of this polynomial can be computed. 
In this way, we can solve the system 
\[
   |q|=1, \quad \Im(Q^3\overline{R}^2)=0,
\]
and prove that a solution $x=\xi$ satisfies the condition
\eqref{eq:swtlQ}
only if $\Re(\xi)=1/2$.
 Substituting $x=\frac{1}{2}+it$ into the second equation
 \eqref{eq:swtlQ}, 
 we have
\begin{equation}\label{eq:non-swallow-fuchs}
   2|R|^4-x(1-x)(2R'Q-RQ')\overline R^2
    =\frac{1}{64}t^2(7-4t^2)^2(21+440t^2-560t^4+256t^6).
\end{equation}
 Since $|q|\neq 1$ at $x=\frac{1}{2}(1\pm\sqrt{7})$, 
 we deduce that the real part of \eqref{eq:non-swallow-fuchs} does not 
 vanish on the singular points.
Hence there is a unique swallowtail in the
image surface of the upper $x$-plane.

\subsubsection{Lambda function}\label{subsub:lambda}
\begin{figure}
\begin{center}
 \includegraphics[width=6cm]{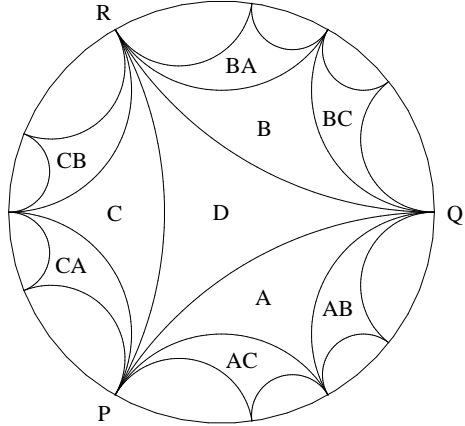}
\end{center}
\caption{Schwarz triangles with three zero angles}
\label{lambda_triangle}
\end{figure}

The inverse of the Schwarz map is a modular function known 
as the lambda function 
\[
  \lambda: \HH^2=\{z\in\C\mid\Im z>0\}\longrightarrow X.
\]
The hyperbolic Schwarz map is expressed in terms of its derivatives. 
In this section we recall its definition and give a few properties.
We begin with the theta functions: for $z\in\HH^2$, 
set $q=e^{\pi iz/2}$,
\[
  \theta_2=\sum_{-\infty}^{\infty}q^{(2n-1)^2/2},\quad
  \theta_3=\sum_{-\infty}^{\infty}q^{2n^2},\quad
  \theta_0=\sum_{-\infty}^{\infty}(-1)^nq^{2n^2}.
\]
Recall the well-known identity
$\theta_3^4-\theta_0^4=\theta_2^4$.
We define the lambda function as
\[
   \lambda(z)=\left(\dfrac{\theta_0}{\theta_3}\right)^4
  =1-16q^2+128q^4-704q^6+3072q^8-11488q^{10}+38400q^{12}-\cdots;
\]
note that $\lambda:\infty\mapsto1,0\mapsto0,1\mapsto\infty$, and 
that $\lambda$ sends every triangle in Figure \ref{lambda_triangle} 
onto the upper/lower half $x$-plane. In the figure, for symmetry reason, 
the Schwarz triangles
tessellating the upper half plane $\HH^2$ are shown in the Poincar\'e disk. 
The inverse of the Schwarz map is given by $x=\lambda(z)$. 
In the expression $\mathscr{S}$ of the hyperbolic Schwarz map 
given in \S \ref{sec:schwarz}, 
the derivatives $\lambda'$ and $\lambda''/\lambda'$ are used. 
They are computed as follows:
Define the Eisenstein series $E_2$ by
\[
  E_2(z)=\frac1{24}\frac{\eta'(z)}{\eta(z)}
  =1-24\sum_{n=1}^\infty\left(\sum_{d|n}d\right)e^{2\pi inz}
  =1-24(q^4+3q^8+4q^{12}+7q^{16}+\cdots).
\]
Then we have 
\[
\dfrac{\theta_0'}{\theta_0}-\dfrac16E_2=-\dfrac16(\theta_2^4+\theta_3^4),\quad
\dfrac{\theta_2'}{\theta_2}-\dfrac16E_2=\dfrac16(\theta_0^4+\theta_3^4),\quad
\dfrac{\theta_3'}{\theta_3}-\dfrac16E_2=-\dfrac16(\theta_0^4-\theta_2^4),
\]
where 
\[
   {}'=q\frac{d}{dq}=\frac2{\pi i}\frac{d}{dz},
\]
and so we have
\[
  \lambda'=-2\theta_2^4\lambda,
\]
which leads to the $q$-series expansion
\[
  \frac{\lambda''}{\lambda'}=(\log \lambda')'=
     4\frac{\theta_2'}{\theta_2}+\frac{\lambda'}{\lambda}
      =\frac46E_2+\frac46(\theta_0^4+\theta_3^4)-2\theta_2^4.
\]
This expression is useful for drawing the picture of the image of $\Phi$, 
because $q$-series converge very fast.

\subsubsection{Visualizing the image surface}
\begin{figure}
\footnotesize
\begin{center}
\begin{tabular}{c@{\hspace{1cm}}c}
 \includegraphics[width=5cm]{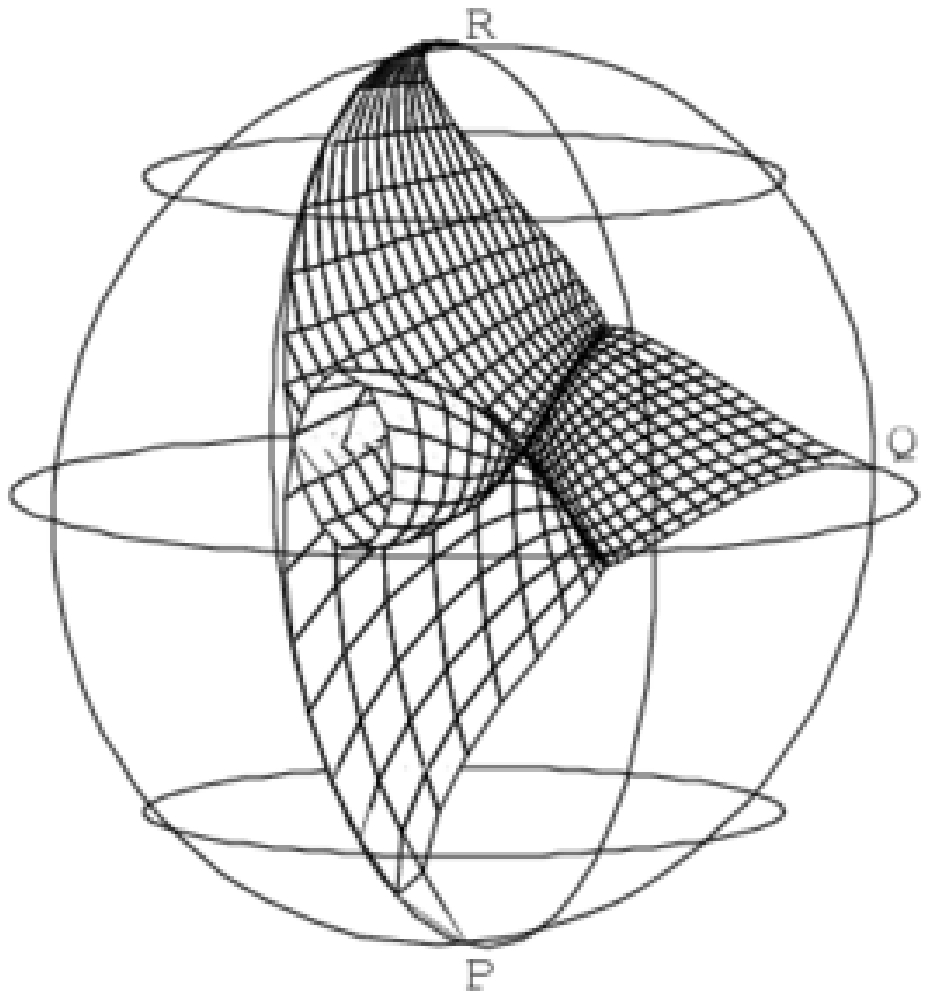} &
 \includegraphics[width=5cm]{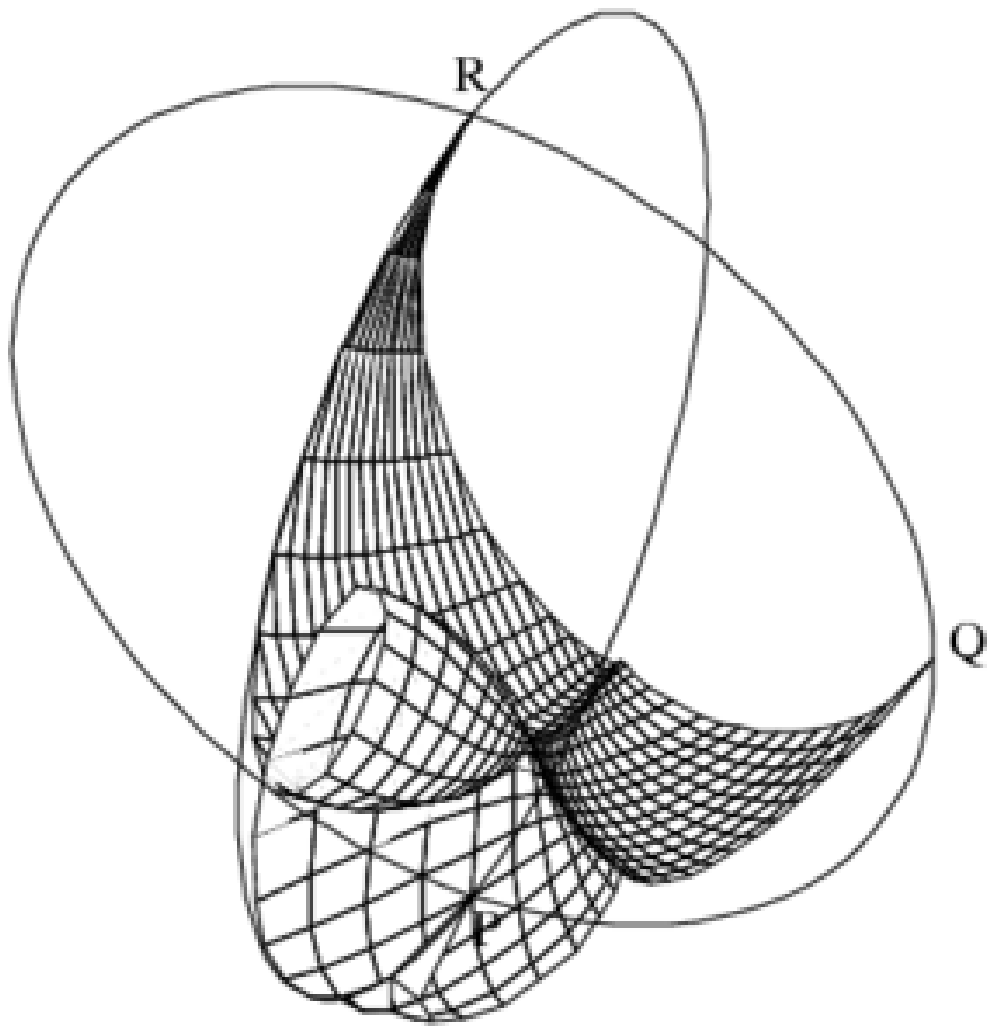} \\
 \multicolumn{2}{c}{Image of one triangle $D$} \\
 \includegraphics[width=5cm]{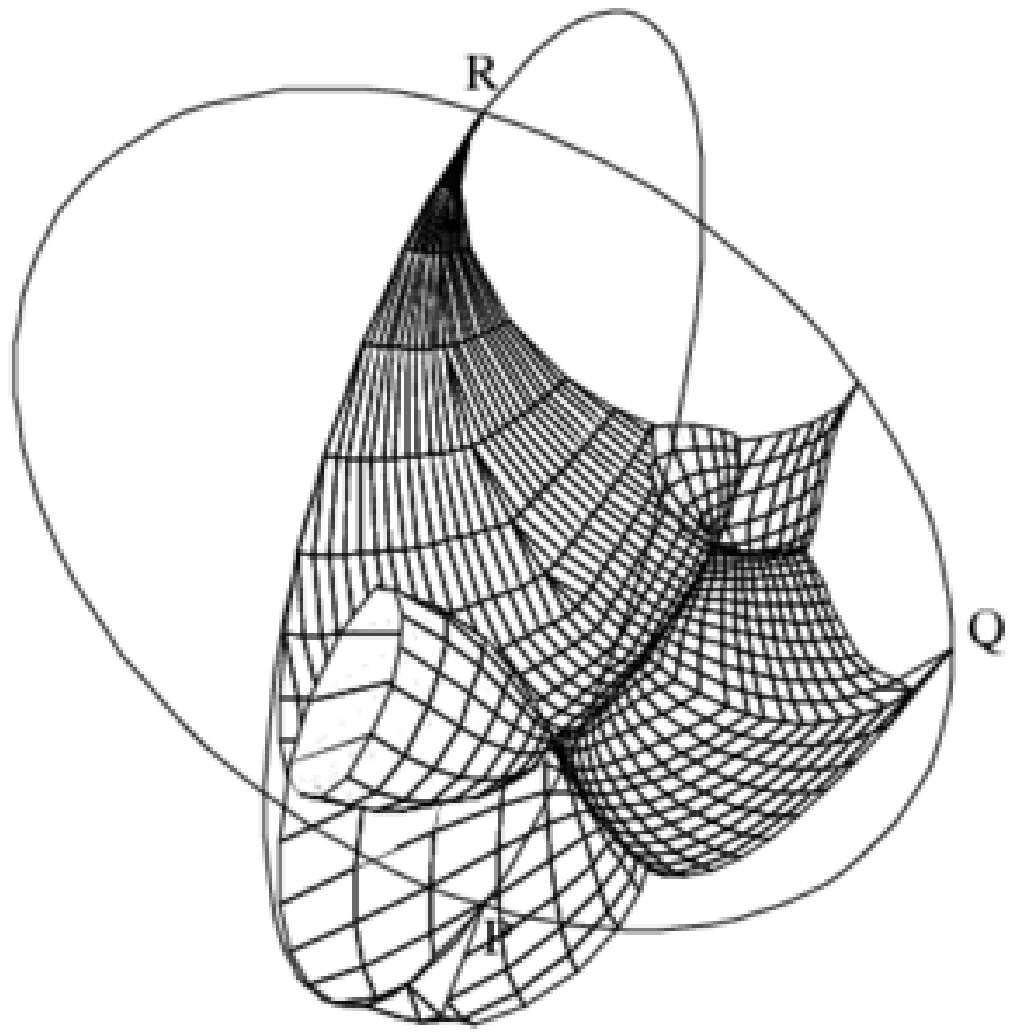} &
 \includegraphics[width=5cm]{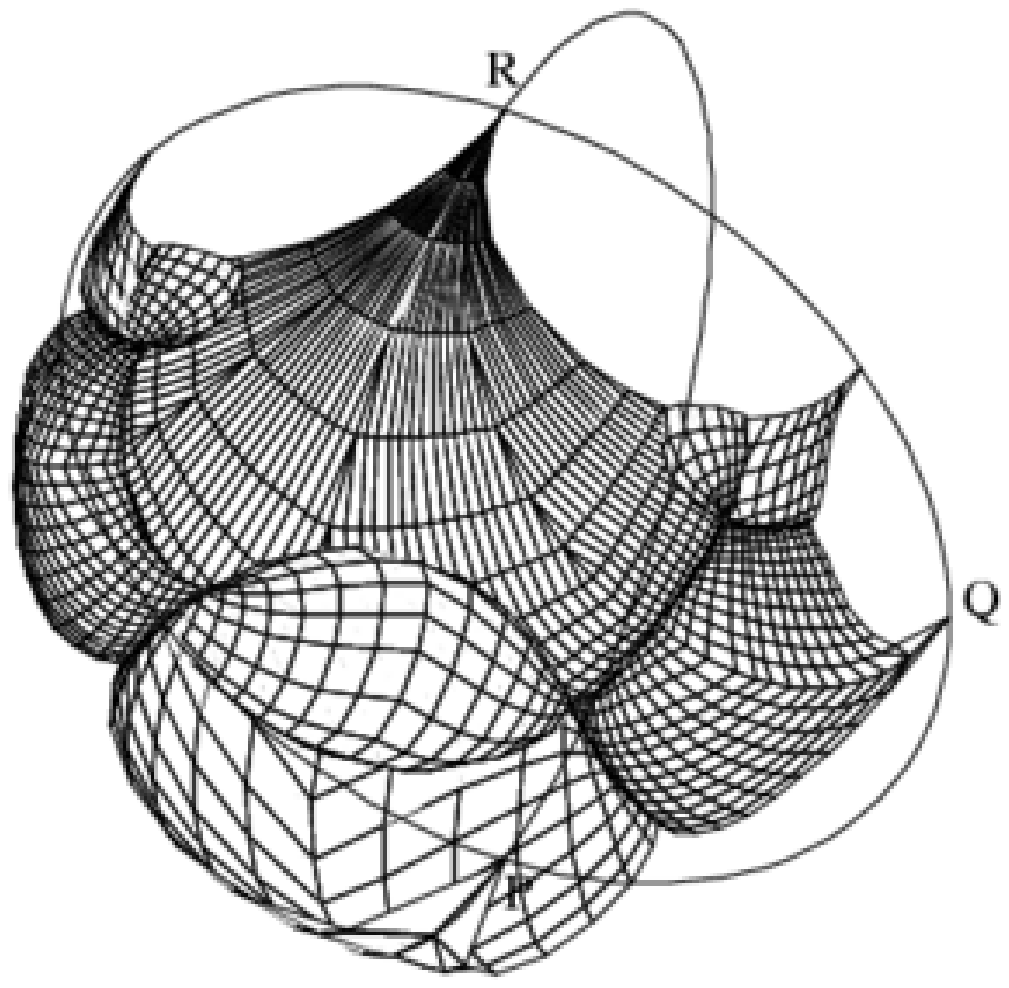} \\
 Image of $\{D,A\}$ &
 Image of $\{D,A,B,BA\}$ \\
 \includegraphics[width=5cm]{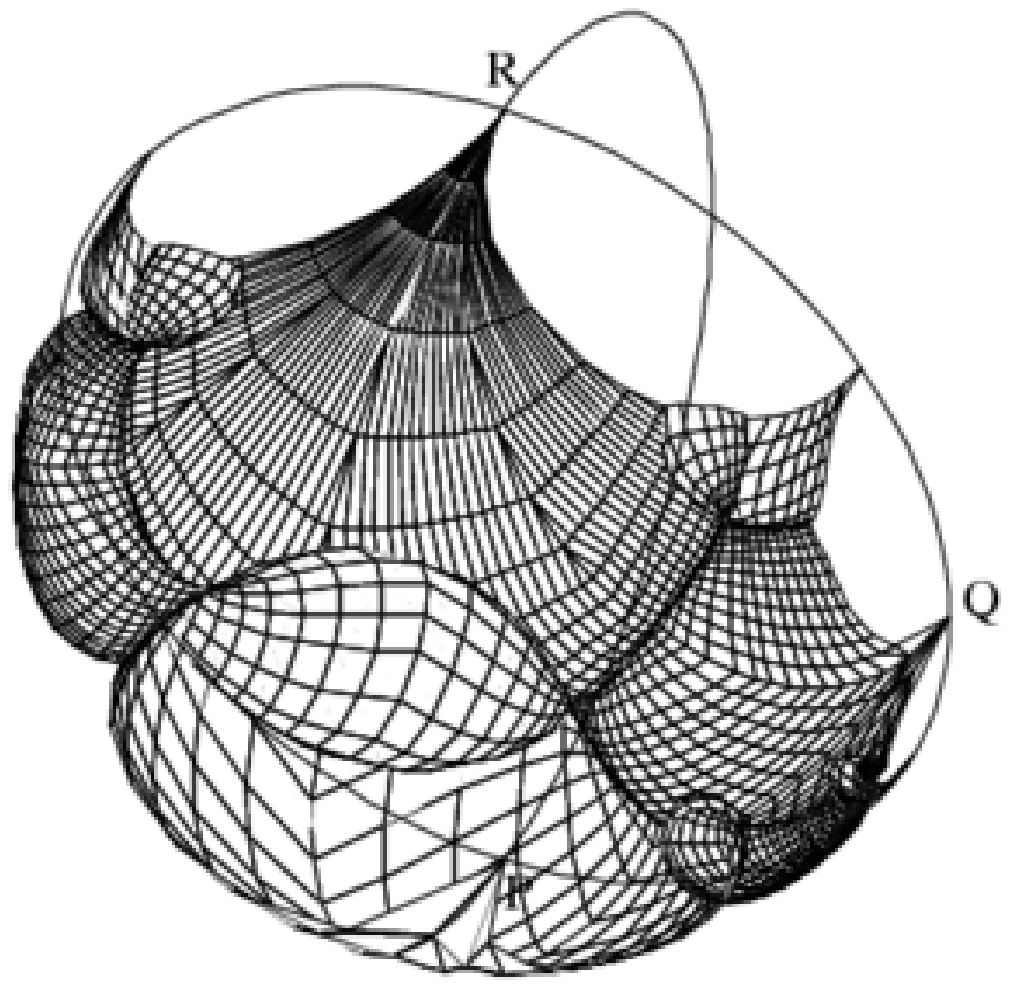} &
 \includegraphics[width=5cm]{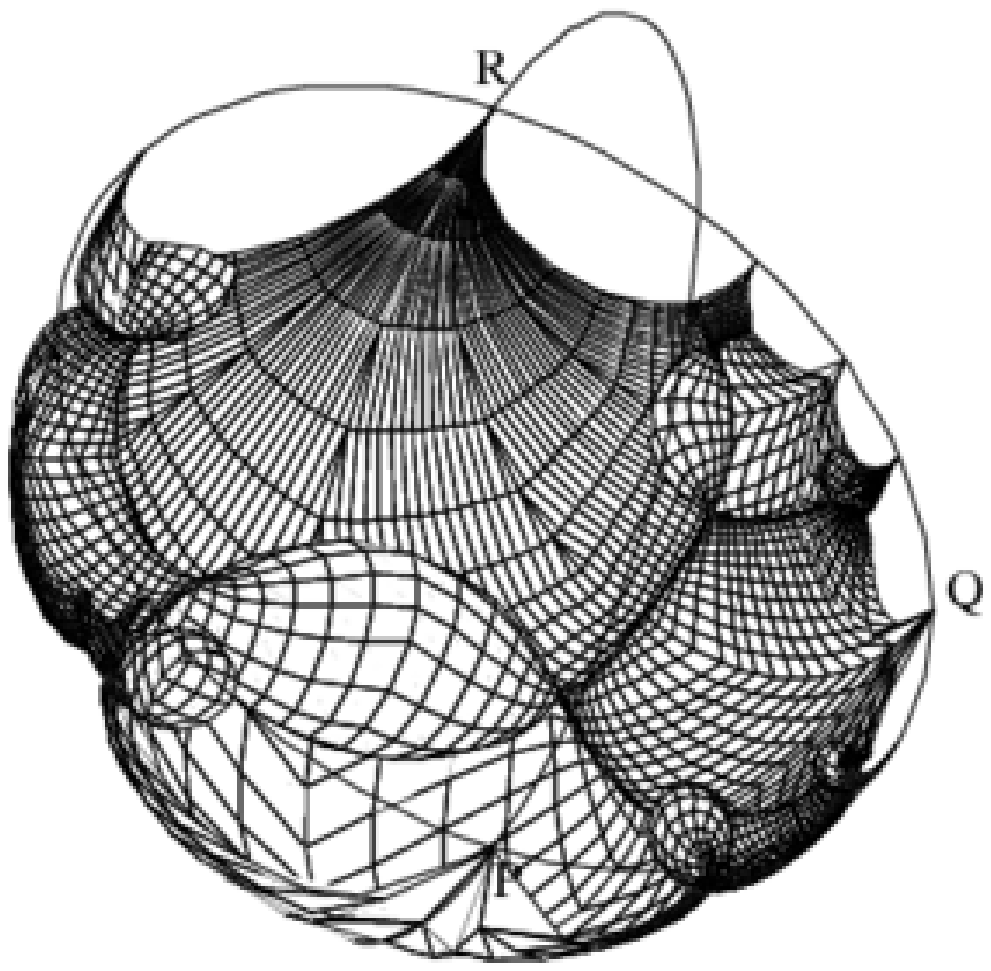} \\
 Image of $\{D,A,B,BA,C,CA\}$ & 
 Image of$\{D,A,B,C,$ \\
 & \multicolumn{1}{r}{$AB,BA,AC,CA,BC,CB\}$}
\end{tabular}
\end{center}
\caption{Images of the hyperbolic Schwarz map\newline 
  when $k_0=k_1=k_\infty=\infty$}
\label{lambda}
\end{figure}

The image of the hyperbolic Schwarz map is shown in 
Figure~\ref{lambda}.
The first one is the image of the triangle $\{D\}$, 
the second one the two triangles $\{D,A\}$, 
the third one the four triangles $\{D,A,B,BA\}$, 
the fourth one the six triangles $\{D,A,B,BA, C, CA\}$, 
and the last one the ten triangles 
$\{D,A,B,C, AB, BA, AC, CA, BC, CB\}$.

\begin{thebibliography}{KR}
\bibitem[Arnold et al. 1985]{AZV}
  V. I. Arnold, S. M. Gusein-Zade, and A. N. Varchenko,
  {\it Singularities of differentiable maps, Vol. $1$},
  Monographs in Math. {\bf 82}, Birkh\"auser, Boston, 1985.
%
\bibitem[Epstein 1986]{Ep} 
  Ch.~L.~Epstein, 
  ``The hyperbolic Gauss map and quasiconformal reflections'',
  {\it J. reine angex. Math.} {\bf 372}(1986), 96--135.
%
\bibitem[G\'alvez et al. 2000]{GMM}
  J. A. G\'alvez, A. Mart\'\i{}nez and F. Mil\'an,
  ``Flat surfaces in hyperbolic $3$-space'',
  {\it Math.\ Annalen} {\bf 316} (2000), 419--435.
%
\bibitem[Iwasaki et al. 1991]{IKSY}
 K. Iwasaki, H. Kimura, S. Shimomura, and M. Yoshida,
  {\it From Gauss to Painlev\'e -- A modern theory of special functions, } 
  Aspects of mathematics {\bf E16}, Vieweg Verlag, Wiesbaden, 1991.
%
\bibitem[Klein 1884]{Kl} F. Klein, 
  {\it Vorlesungen \"uber das Ikosaeder
  und die aufl\"osung der Gleichungen vom f\"unften Grade, }
  Teubner, Leipzig, 1884.
%
\bibitem[Kokubu et al. 2003]{KUY} 
  M. Kokubu, M. Umehara, and K. Yamada,
  ``An elementary proof of Small's formula for null curves in $PSL(2,\C)$ 
  and an analogue for Legendreian curves in $PSL(2,\C)$'', 
  Osaka J. Math. {\bf 40}(2003), 697--715.
%
\bibitem[Kokubu et al. 2004]{KUY2} 
  M. Kokubu, M. Umehara, and K. Yamada,
  ``Flat fronts in hyperbolic $3$-space'',
  Pacific J. Math. {\bf 216}(2004), 149--175.
%
\bibitem[Kokubu et al. 2005]{KRSUY} 
  M. Kokubu, W. Rossman, K. Saji, M. Umehara, and K. Yamada,
  ``Singularities of flat fronts in  hyperbolic space'', 
  Pacific J. Math. {\bf 221}(2005), 303--351.
%
\bibitem[Noro et al.  2007]{NSYY} 
  M. Noro, T. Sasaki, K. Yamada and M. Yoshida,
``Confluence of swallowtail singularites of the hyperbolic
Schwarz map defined by the hypergeometric differential equation'',
preprint, 2007.
%
\bibitem[Sasaki et al.  2007]{SYY2} 
  T. Sasaki, K. Yamada and M. Yoshida,
  ``Derived Schwarz map
  of the hypergeometric differential equation 
  and a parallel family of flat fronts'',
  preprint, 2007.
%
\bibitem[Yoshida 1997]{Yo}
  M. Yoshida,
  {\it Hypergeometric Functions, My Love,} 
  Aspects of mathematics {\bf E32}, Vieweg Verlag, Wiesbaden, 1997.
\end{thebibliography}
\end{document}